\documentclass[12pt]{article}
\usepackage[numbers,sort&compress]{natbib}
\usepackage{amsmath,amssymb}
\usepackage{graphicx}
\usepackage{amsthm}



\numberwithin{equation}{section}

\newcommand{\RA}{\rightarrow}
\newcommand{\LS}{\limsup_{t\rightarrow\infty}}
\newcommand{\LI}{\liminf_{t\rightarrow\infty}}
\newcommand{\D}{{\rm d}}
\newtheorem{definition}{\bf Definition}[section]
\newtheorem{corollary}{\bf Corollary}[section]
\newtheorem{remark}{\bf Remark}[section]

\newtheorem{lemma}{\bf Lemma}[section]
 \newtheorem{theorem}{\bf Theorem}[section]
\begin{document}
\title{\bf Some new results on the dynamical behavior of non-autonomous logistic system with
random perturbation}
\author{Hongxiao Hu\footnote{Corresponding author. Tel.: +86 15000559112.
E-mail address: hhxiao1@126.com (H.X. Hu).} \\
\small College of Science,  University of Shanghai for Science and Technology, \\
\small Shanghai 200093, P. R. China}
\date{}
\maketitle
\baselineskip 18pt
\textbf{Abstract.} In this paper, a non-autonomous stochastic logistic system is considered. An interesting result on the effect of stochastically perturbation for the dynamic behavior are obtained. That is, under certain conditions the stochastic system have similar dynamic behave with the deterministic system, but if the noise is sufficiently large, it will spoil these nice properties.
Furthermore, we introduce a new research method for studying the stochastic equation and some new sufficient conditions for the stochastic bounded, stochastic permanence, extinction and global attractivity of the system are established.
The corresponding results (see, for example, references [DQ Jiang, NZ Shi, and XY Li. Global stability and stochastic permanence of a
non-autonomous logistic equation with random perturbation. J Math Anal Appl,
340:588-597, 2008; M Liu and K Wang. Persistence and extinction in stochastic non-autonomous logistic
systems. J Math Anal Appl, 375:443-457, 2011; XY Li and XR Mao. Population dynamical behavior of non-autonomous lotka-
volterra competitive system with random perturbation.
Discrete Contin Dyn
Syst, 24:523-545, 2009]) are improved and extended.


\section{Introduction}
It is an usual phenomenal in nature that a single species whose members compete among themselves for limited resources. Therefore, the investigation of the single species model is one of the dominant themes in both ecology and mathematical ecology. As we know, one of the famous models for single species system is the classical non-autonomous logistic system which can be expressed as follows:
\begin{equation}
\frac{\D x(t)}{\D t}=x(t)[(r(t)-a(t)x(t))],\label{S2}
\end{equation}
where $x(t)$ denotes the population size of the species at time $t$, $r(t)$ and $a(t)$ are continuous bounded functions on $[0,\infty)$, and $a(t)$ is nonnegative. We refer the reader to \cite{May} for a detailed model construction. Owing to its theoretical and practical significance, system (\ref{S2}) and its generalization form have received great attention and has been studied extensively (see \cite{Teng21,FW,KG,YK,TL,Teng,BL}).

On the other hand, many species are always affected inevitably by environmental noise which is an important component in an ecosystem (see \cite{G84,G86}). In resent years, many authors consider the effect of random disturbance in population dynamics (see \cite{M4,G84,G86,M20,M25,M26, M5,M23,M6,M12,JSL,LM,LW} and the reference therein). In particular, \cite{M23} showed that different structures of environmental noise may have different effects on the population system. Currently we have two kinds of systems to model the effect of the environmental fluctuations in population dynamics. one is the white noise affects $a(t)$ mainly, i.e., the intraspecific competition coefficient is stochastically perturbed (see \cite{M25,M26,M5,W12,W13,LW}). The other is the growth rate $r(t)$ is subject to environmental noise (see \cite{W14,W15,W16,JSL,W18,W19,LM,W22,W23,LW}). If we assume that the growth rate $r(t)$ is stochastically perturbed with
$$
    r(t)\RA r(t)+\sigma(t)\dot{B}(t),
$$
where $\dot{B}(t)$ is a white noise and $\sigma^2(t)$ represents the intensity of the noise. Then the stochastically perturbed system can be described by the following It$\rm\hat{o}$ equation:
\begin{equation}
\D x(t)=x(t)[(r(t)-a(t)x(t))\D t+\sigma(t)\D B(t)],\label{S1}
\end{equation}
where $B(t)$ is standard Brownian motions, $r(t)$, $a(t)$ and $\sigma(t)$ are all continuous bounded functions on $[0,+\infty)$, $a(t)$ and $\sigma(t)$ are nonnegative.

In \cite{JSL}, the authors show that system (\ref{S1}) is stochastically permanent and globally attractive provided $r(t)$, $a(t)$ and $\sigma(t)$ are continuous $T$-periodic functions, $r(t)>0$, $a(t)>0$ and $\min_{t\in[0,T]}r(t)>\max_{t\in[0,T]}\sigma^2(t)$. \cite{LM} improved the results in \cite{JSL} and obtained that if $\min_{t\in[0,\infty)}a(t)>0$ and $\min_{t\in[0,\infty)}\{r(t)-\frac{1}{2}\sigma^2(t)\}>0$, then system (\ref{S1}) is permanence. In \cite{LW}, stochastically permanence of system (\ref{S1}) is obtained with $a(t)>0$ and $\LI (r(t)-\frac{1}{2}\sigma^2(t))>0$. An indication of the generality of these results emerges from biological interpretations of the assumptions. The principal assumptions of these results are that the growth rate must be greater than $\frac{1}{2}\sigma^2(t)$ and the carrying capacity must be positive all the time. But as we known, the growth properties of every natural population is affected by the population's environment. Physical environmental conditions usually change greatly through the year and can influence species directly. Good weather can stimulate growth in body size and bad weather can cause death. Hence, it is difficult to keep the growth rate and the carrying capacity positive in all time. But we can find some fixed length of time such that during any time interval of this length the growth rate and carrying capacity are positive. In \cite{TL}, the authors studied system (\ref{S2}). They obtained permanence and global asymptotic stability of system (\ref{S2}) with the conditions that the growth rate and carrying capacity are positive during any time interval of fixed length after some fixed time. However, we see that, up to now, the same results still have not been established for the stochastic non-autonomous logistic equation.

Motivated by above works, in this paper, our purpose is namely to extend the results given in \cite{TL} to the stochastic non-autonomous logistic system (\ref{S1}). We will establish a series of very general and rather weak criteria on the stochastic ultimately bounded, stochastically permanence, extinction and global asymptotic stability for system (\ref{S1}). We will see that when system (\ref{S1}) degenerate into the deterministic system (\ref{S2}) these criteria will be very similar to the corresponding results given in \cite{TL}. And we will also see that the corresponding results obtained by \cite{JSL,LM,LW} are improved and extended in these criteria. Therefore, the results obtained in this paper are completely new and useful.

\section{Preliminaries}
Let $(\Omega,\mathcal{F},P)$ be a complete probability space with a filtration $\{\mathcal{F}_t\}_{t\ge0}$ satisfying the usual conditions, that is it is right continuous and increasing with $\mathcal{F}_0$ contains all $P$-null sets. Let $B(t)$ denote the standard Browian motions defined on this probability. We also denote the positive number by $R_+$, that is
$$
    R_+=\{x\in R:x>0\}.
$$
If $f(t)$ is a continuous bounded function on $[0,+\infty)$, we denote
$$
    f_u=\sup_{t\in[0,+\infty)}f(t)\quad\mbox{and}\quad f_l=\inf_{t\in[0,+\infty)}f(t).
$$
For any set $A$, we denote $I_A$ be the indicator function of $A$.
For system (\ref{S1}) we introduce
the following assumptions:

($\rm H_1$) There is a  positive constant $\gamma$ such that
    $$
        \LI\int_t^{t+\gamma}a(s)\D s>0.
   $$

($\rm H_2$) There is a  positive constant $\lambda$ such that
    $$
       \LI\int_t^{t+\lambda}\big(r(s)-\frac{1}{2}\sigma^2(s)\big)\D s>0.
    $$

($\rm H_3$) There is a  positive constant $\lambda$ such that
    $$
       \LS\int_t^{t+\lambda}\big(r(s)-\frac{1}{2}\sigma^2(s)\big)\D s\le0.
    $$
\begin{remark}
    \rm If $a_l$ and $(r-\frac{1}{2}\sigma^2)_l$ are positive, then ($\rm H_1$) and ($\rm H_2$) are satisfied. And ($\rm H_3$) hold if $(r-\frac{1}{2}\sigma^2)_l$ is negative. Therefore, the three assumptions are very general and weak.
\end{remark}

It is well known that, in order for a stochastic differential equation to have an unique global solution for any given initial data, the coefficients of this equation are generally required to satisfy linear growth condition and the local Lipschitz condition. Obviously, the coefficients of equation (\ref{S1}) do not satisfy the linear growth condition. From the Theorem 2.1 in \cite{LM}, we can obtain that equation (\ref{S1}) has a global positive solution, that is
\begin{lemma}\label{L1}
    For any given initial value $x_0\in R_+$, there is an unique solution $x(t)$ to equation (\ref{S1}) on $t\ge0$ and the solution will remain in $R_+$ with probability one.
\end{lemma}
In the following, we give out a useful lemma which will be used in the next section.
\begin{lemma}\label{L2}
    Let $X=\{X_t;0\le t<\infty\}$ is a stochastic process on the probability space $(\Omega,\mathcal{F},P)$, for any positive constance $m$ and $M$ we have

    (i) $\LS P\{X_t>M\}\le P\{\LS X_t>M\}$;

    (ii) $\LS P\{X_t<m\}\le P\{\LI X_t<m\}$.\\
\rm\textbf{Proof.} From the well known Fatou Lemma we have
\begin{equation}
\LS P\{X_t>M\}=\LS\int_\Omega I_{\{X_t>M\}}\D P\le\int_\Omega\LS I_{\{X_t>M\}}\D P\label{L2.1}.
\end{equation}
For any $\omega\in\Omega$, if $\LS I_{\{X_t>M\}}(\omega)=1$, then by the definition of superior limit and indicator function there is a time sequence $\{t_k\}$ with $t_k\RA\infty$ as $k\RA\infty$ such that
$$
    I_{\{X_{t_k}>M\}}(\omega)=1\quad\mbox{for all } k=1,2,\cdots.
$$
This means
$$
    \omega\in\{X_{t_k}>M\}\quad\mbox{for all } k=1,2,\cdots.
$$
Hence,
$$
    \omega\in\{\LS X_t>M\}.
$$
Therefore,
$$
\int_\Omega\LS I_{\{X_t>M\}}\D P\le P\{\LS X_t>M\}.
$$
Together with inequality (\ref{L2.1}), (i) hold obviously.

In the following we will prove the case (ii) of this lemma.
By using the Fatou Lemma, we also have
\begin{equation}
\LS P\{X_t<m\}=\LS\int_\Omega I_{\{X_t<m\}}\D P\le\int_\Omega\LS I_{\{X_t<m\}}\D P\label{L2.2}.
\end{equation}
For any $\omega\in\Omega$, if $\LS I_{\{X_t<m\}}(\omega)=1$, by using the similar argument as the case (i) we have a time sequence $\{t_k\}$ with $t_k\RA\infty$ as $k\RA\infty$ such that
$$
    \omega\in\{X_{t_k}<m\}\quad\mbox{for all } k=1,2,\cdots.
$$
Hence,
$$
    \omega\in\{\LI X_t<m\}.
$$
Therefore,
$$
\int_\Omega\LS I_{\{X_t<m\}}\D P\le P\{\LI X_t<m\}.
$$
It follows equation (\ref{L2.2}) that (ii) is true. This complete the proof.\qed
\end{lemma}
\begin{remark}
    \rm Since
    $$\LS P\{X_t>M\}=1-\LI P\{X_t\le M\}$$
    and
    $$P\{\LS X_t>M\}=1-P\{\LS X_t\le M\},$$
    from (i) of above lemma we can obtain
    $$
        \LI P\{X_t\le M\}\ge P\{\LS X_t\le M\}.
    $$
    Similarly, from (ii) of above lemma we can have
    $$
        \LI P\{X_t\ge m\}\ge P\{\LI X_t\ge m\}.
    $$
\end{remark}
Let $M(t)=\int_0^t\sigma(s)\D B(s)$ is a martingale, the quadratic variation of this martingale is
$$
    \langle M\rangle_t=\int_0^t\sigma^2(s)\D s\le\sigma^2_ut.
$$
By the strong law of large numbers or martingale (see \cite{M97}), we therefore have
\begin{equation}
    \lim_{t\RA\infty}\frac{M(t)}{t}=0\quad a.s.\label{P1}
\end{equation}
Let
$$
    \Omega_0=\{\lim_{t\RA\infty}\frac{M(t)}{t}=0\},
$$
obviously, $P(\Omega_0)=1$.
\section{Stochastically permanence}
In this section, we will study the stochastic ultimate boundedness and stochastically permanence which are defined as follows.
\begin{definition}
    \rm Equation (\ref{S1}) is said to be stochastically ultimately bounded, if for any $\varepsilon\in(0,1)$, there exists a positive constant $M=M(\varepsilon)$ such that
    $$
       \LI P\{x(t)\le M\}\ge1-\varepsilon,
    $$
    for any positive solution $x(t)$ of system (\ref{S1}).
\end{definition}
\begin{definition}
    \rm Equation (\ref{S1}) is said to be stochastically permanent if for any $\varepsilon\in(0,1)$, there exists a pair of positive constants $m=m(\varepsilon)$ and $M=M(\varepsilon)$ such that
    $$
        \LI P\{x(t)\ge m\}\ge1-\varepsilon\mbox{ and } \LI P\{x(t)\le M\}\ge1-\varepsilon,
    $$
    for any positive solution $x(t)$ of system (\ref{S1}).
\end{definition}

 Firstly, we have the following theorem on the asymptotically bounded in any $p$th moment of the solutions of equation (\ref{S1}).
\begin{theorem}\label{B}
    Suppose ($H_1$) hold, then for any positive constant $p$ there exist positive constant $L(p)$ such that
    $$
        \LS E x^p(t)\le L(p)
    $$
    for any positive solution $x(t)$ of equation (\ref{S1}).\\
    \rm\textbf{Proof.} By using the similar argument as the Lemma 2.3 in \cite{JSL}, we can obtain
    $$
        \frac{\D E[x^p(t)]}{\D t}\le pE[x^p(t)]\Big(\big(r(t)+\frac{1}{2}(p-1)\sigma^2(t)\big)-a(t)\big(E[x^p(t)]\big)^\frac{1}{p}\Big)
    $$
    Let $y(t)=\big(E[x^p(t)]\big)^\frac{1}{p}$, then we have
    $$
        \frac{\D y(t)}{\D t}\le y(t)\Big(\big(r(t)+\frac{1}{2}(p-1)\sigma^2(t)\big)-a(t)y(t)\Big).
    $$
    We consider the following auxiliary equation
     \begin{equation}
        \frac{\D z(t)}{\D t}=z(t)\Big(\big(r(t)+\frac{1}{2}(p-1)\sigma^2(t)\big)-a(t)z(t)\Big).\label{A1}
     \end{equation}
    By the assumption ($\rm H_1$) and Lemma 1 in \cite{TL}, we obtain that there is a positive constant $M$ such that
    $$
        \LS z(t)\le M
    $$
    for any positive solution $z(t)$ of equation (\ref{A1}). By the comparison theorem, we have
    $$
        \LS y(t)\le \LS z(t)\le M
    $$
    with initial value $y(0)=z(0)$. Therefore,
    $$
        \LS E[x^p(t)]=\LS y^p(t)\le M^p:=L(p).
    $$
    This complete the proof.\qed
\end{theorem}
By Chebyshev's inequality and Theorem \ref{B}, the following result is straightforward.
\begin{theorem}\label{B2}
    suppose ($H_1$) hold, then system (\ref{S1}) is stochastically ultimately bounded.
\end{theorem}
\begin{remark}
    \rm In \cite{JSL,LM}, the authors obtained the stochastically ultimately bounded of system (\ref{S1}) with $a_l>0$. Obviously the conditions of Theorem \ref{B2} is more weaker than these. Furthermore, we can obtain more general result under the condition ($\rm H_1$), which we will discuss in the following theorem.
\end{remark}
\begin{theorem}\label{T1}
    Suppose ($H_1$) hold, then for any $\varepsilon\in(0,1)$ there is a positive constant $M=M(\varepsilon)$ such that for any initial value $x_0\in R_+$ the solution obeys
    $$
        P\{\LS x(t)\le M\}\ge 1-\varepsilon.
    $$
\rm\textbf{Proof.} From the assumption ($\rm H_1$), there are positive constants $T_0$, $L>1$ and $\mu$ such that
\begin{equation}
    \int_t^{t+\gamma}\big(r(s)-\frac{1}{2}\sigma^2(s)-a(s)L\big)\D s<-\mu\quad\mbox{for all }t\ge T_0.\label{T1.1}
\end{equation}
Then we can choose a positive constant $\varepsilon_0\ll1$ such that
\begin{equation}
    \int_t^{t+\gamma}\big(r(s)-\frac{1}{2}\sigma^2(s)-a(s)(L-\varepsilon_0)\big)\D s<-\frac{\mu}{2}\quad\mbox{for all }t\ge T_0.\label{T1.1'}
\end{equation}
By the It$\rm\hat{o}$ formula, we have
$$
    \D \ln x(t)=[r(t)-\frac{1}{2}\sigma^2(t)-a(t)x(t)]\D t+\sigma(t)\D B(t).
$$
For any two positive constant $t_1<t_2$, integrating above equation from $t_1$ to $t_2$ we have
\begin{equation}
    \ln x(t_2)-\ln x(t_1)=\int_{t_1}^{t_2}[r(t)-\frac{1}{2}\sigma^2(t)-a(t)x(t)]\D t+\int_{t_1}^{t_2}\sigma(t)\D B(t).\label{T1.4}
\end{equation}
Firstly, we will prove that
\begin{equation}
    \LI x(t,x_0,\omega)<L\quad\mbox{for all }x_0\in R_+\mbox{ and }\omega\in\Omega_0.\label{CL1}
\end{equation}
If it is not true, then there exist a $x_0\in R_+$ and $\omega_0\in\Omega_0$ such that
$$
    \LI x(t,x_0,\omega_0)\ge L.
$$
It follows from this that for above $\varepsilon_0$ there is a $T_1\ge T_0$ such that
\begin{equation}
    x(t,x_0,\omega_0)\ge L-\varepsilon_0\quad\mbox{ for all }t\ge T_1.\label{T1.5}
\end{equation}
By (\ref{T1.4}), we have
$$
    \ln x(t)-\ln x(T_1)=\int_{T_1}^t\big(r(s)-\frac{1}{2}\sigma^2(s)-a(s)x(s)\big)\D s+M(t)-M(T_1).
$$
Choosing a positive integer $p$ such that $t\in[T_1+p\gamma,T_1+(p+1)\gamma)$, then by (\ref{T1.1'}) and (\ref{T1.5}) we get
\begin{eqnarray}
    \ln x(t,\omega_0)-\ln x(T_1,\omega_0)&\le&\big(\int_{T_1}^{T_1+p\gamma}+\int_{T_1+p\gamma}^t\big)\big(r(s)-\frac{1}{2}\sigma^2(s)-a(s)(L-\varepsilon_0)\big)\D s\nonumber\\
    &&+M(t,\omega_0)-M(T_1,\omega_0)\nonumber\\
    &\le&-\frac{1}{2}p\mu+\gamma\alpha_1+M(t,\omega_0)-M(T_1,\omega_0)\nonumber\\
    &\le&-\frac{\mu}{2\gamma}(t-T_1-\gamma)+\gamma\alpha_1+M(t,\omega_0)-M(T_1,\omega_0),\label{T1.6}
\end{eqnarray}
where $\alpha_1=\sup_{s\ge0}|r(s)-\frac{1}{2}\sigma^2(s)-a(s)(L-\varepsilon_0)|$. Dividing $t$ on the both sides and then letting $t\RA\infty$, it finally follows from (\ref{P1}) that
$$
    \LS\frac{\ln x(t)}{t}\le-\frac{\mu}{2\gamma}.
$$
This implies $\lim_{t\RA\infty}x(t)=0$, which is a contradiction. Therefore, (\ref{CL1}) hold.

In the following, we will prove this theorem is true. Otherwise, there is a positive constant $\delta_0\in(0,1)$ such that for any positive integer $n\ge2$ we have a $x_n\in R_+$ satisfying the property that
$$
    P\{\LS x_n(t)\le nL\}\le1-\delta_0,
$$
where $x_n(t)=x(t,x_n)$. Consequently,
\begin{equation}
    P\{\LS x_n(t)>nL\}\ge \delta_0\label{T1.7}
\end{equation}
Let us now define a sequence of stopping time for any $n\ge2$,
\begin{eqnarray*}
    \sigma_1^{(n)}&=&\inf\{t\ge T_0:x_n(t)\ge nL\},\quad \sigma_{2k}^{(n)}=\inf\{t\ge \sigma_{2k-1}^{(n)}:x_n(t)\le L\}\\
    \sigma_{2k+1}^{(n)}&=&\inf\{t\ge \sigma_{2k}^{(n)}:x_n(t)\ge nL\},\quad\mbox{for }k=1,2,\cdots.
\end{eqnarray*}
Let
$$\Omega_n=\Big(\{\LS x_n(t)> nL\}\cup\Big(\{\LS x_n(t)=nL\}\cap\big(\bigcap_{k=1}^{\infty}\{\sigma_k^{(n)}<\infty\}\big)\Big)\Big)\cap\Omega_0$$ for $n\ge2.$
 It follows from (\ref{T1.7}) that
\begin{equation}
    P(\Omega_n)\ge\delta_0.\label{T1.8}
\end{equation}
Note from (\ref{CL1}) and the definition of $\Omega_n$ that for any $n=2,3,\cdots$
$$
    \sigma_k^{(n)}<\infty\quad\mbox{for any }k\ge 1\mbox{ whenever }\omega\in\Omega_n,
$$
$\{\sigma_k^{(n)}<\infty\}$ is non-increasing with $k$ and
\begin{equation}
    \Omega_n=\Big(\bigcap_{k=1}^{\infty}\{\sigma_k^{(n)}<\infty\}\Big)\cap\Omega_0.\label{T1.9}
\end{equation}
Now for any $n\ge2$ we can define a time sequence $\{\tau_k^{(n)}\}$, where
$$
    \tau_k^{(n)}=\sup\{t:x_n(t)\le L\mbox{ and }\sigma_{2k}^{(n)}\le t<\sigma_{2k+1}^{(n)}\},\quad k=1,2,\cdots.
$$
Since $x_n(t)$ is adapted to $\{\mathcal{F}_t\}$ for any $n$, we have that $\tau_k^{(n)}$ is $\mathcal{F}$-measurable but is not a stopping time for any $k$ and $n$. By (\ref{CL1}) and the definition of $\Omega_n$, we can find that for any $\omega\in\Omega_n$
\begin{equation}
    \begin{array}{c}
        L<x_n(t)<nL\quad\mbox{for all }t\in(\tau_k^{(n)},\sigma_{2k+1}^{(n)}),\vspace{1mm}\\

      x_n(\tau_k^{(n)})=L \mbox{ and }x_n(\sigma_{2k+1}^{(n)})=nL\quad\mbox{for all }n\ge2,\;k=1,2,\cdots.
    \end{array}\label{T1.11}
\end{equation}

Let $\alpha_2=\sup_{s\ge0}|r(s)-\frac{1}{2}\sigma^2(s)-a(s)L|$ and $\Omega_k^n(T)=\{\sigma_{2k+1}^{(n)}-\tau_k^{(n)}\ge T\}\cap\Omega_n$ for all $k$ and $n$. We claim that for any positive constant $T>\frac{\gamma}{\mu}(\mu+\alpha_2\gamma+1)$, there are integers $N$ and $L_n$ such that
\begin{equation}
    P\left\{\Omega_k^n(T)\right\}=P(\Omega_n)\quad\mbox{for all }n\ge N, k\ge L_n.\label{CL2}
\end{equation}
If the claim is not true, there is a positive constant $T_2>\frac{\gamma}{\mu}(\mu+\alpha_2\gamma+1)$, for any $i=1,2,\cdots$ there exist $n_i>i$ such that for any $h=1,2,\cdots$ we have a $k_i^h\ge h$ satisfying
\begin{equation}
    P\big\{\Omega_{k_i^h}^{n_i}(T_2)\big\}<P(\Omega_{n_i}).\label{TT}
\end{equation}
By (\ref{T1.9}), for any positive constant $\rho<\delta_0^2/(9\sigma_u^2T_2)$ and $n\ge2$ there exist a $K_n=K_n(\rho)$ such that
\begin{equation}
    P\big\{\big(\{\sigma_k^{(n)}<\infty\}\cap\Omega_0\big)\backslash\Omega_n\big\}<\rho\quad\mbox{for all }k\ge K_n.\label{T1.10}
\end{equation}
For $n_i$ we can choose $h=K_{n_i}$ then there is a $k_i:=k_i^h\ge h$ satisfying (\ref{TT}) and (\ref{T1.10}).
Consequently, it follows from (\ref{TT}) we have
$$
    P\big\{\{\sigma_{2k_i+1}^{(n_i)}-\tau_{k_i}^{(n_i)}<T_2\}\cap\Omega_{n_i}\big\}>0.
$$
Obviously, there are two cases:\\
Case 1. There is a positive constant $\eta_0$ such that
$$
    P\big\{\{\sigma_{2k_i+1}^{(n_i)}-\tau_{k_i}^{(n_i)}<T_2\}\cap\Omega_{n_i}\big\}>\eta_0\quad\mbox{for all }i=1,2,\cdots;
$$
Case 2. There are subsequences $\{n_{i_j}\}\subset\{n_i\}$ and $\{k_{i_j}\}\subset\{k_i\}$, a sequence $\{\eta_j\}$ with $\eta_j\RA 0$ as $j\RA\infty$ and $\eta_j<\delta_0$ for all $j=1,2, \cdots$ such that
\begin{equation}
    0<P\big\{\{\sigma_{2k_j+1}^{(n_j)}-\tau_{k_j}^{(n_j)}<T_2\}\cap\Omega_{n_j}\big\}<\eta_j\quad\mbox{for all }j=1,2,\cdots,\label{TC2}
\end{equation}
where $k_j:=k_{i_j}$ and $n_j:=n_{i_j}$.

In the following, we will prove neither Case 1 or Case 2 is true. If Case 1 arises,
 follows from (\ref{T1.4}), (\ref{T1.11}) and the moment inequality of stochastic integrals (see \cite{M97}), we compute
\begin{eqnarray*}
    \eta_0\ln n_i&=&\eta_0(\ln n_iL-\ln L)\\
    &\le&E\big[\sup_{0\le t\le T_2}I_{\{\sigma_{2k_i+1}^{(n_i)}-\tau_{k_i}^{(n_i)}<T_2\}\cap\Omega_{n_i}}\big(\ln x_{n_i}(\sigma_{2k_i+1}^{(n_i)})-\ln x_{n_i}(\sigma_{2k_i+1}^{(n_i)}-t)\big)\big]\\
    &\le&E\big[\sup_{0\le t\le T_2}I_{\{\sigma_{2k_i+1}^{(n_i)}-\tau_{k_i}^{(n_i)}<T_2\}\cap\Omega_{n_i}}
    \int_{\sigma_{2k_i+1}^{(n_i)}-t}^{\sigma_{2k_i+1}^{(n_i)}}\big(r(s)-\frac{1}{2}\sigma^2(s)-a(s)x(s)\big)\D s\big]\\
    &&+E\big[\sup_{0\le t\le T_2}I_{\{\sigma_{2k_i+1}^{(n_i)}-\tau_{k_i}^{(n_i)}<T_2\}\cap\Omega_{n_i}}
    \int_{\sigma_{2k_i+1}^{(n_i)}-t}^{\sigma_{2k_i+1}^{(n_i)}}\sigma(s)\D B(s)\big]\\
    &\le&r_uT_2P\big(\{\sigma_{2k_i+1}^{(n_i)}-\tau_{k_i}^{(n_i)}<T_2\}\cap\Omega_{n_i}\big)\\
    &&+E\Big[I_{\{\sigma_{2k_i+1}^{(n_i)}<\infty\}}\sup_{0\le t\le T_2}\big|\int_{\sigma_{2k_i+1}^{(n_i)}-t}^{\sigma_{2k_i+1}^{(n_i)}}\sigma(s)\D B(s)\big|\Big]\\
    &\le&r_uT_2+8E\Big[I_{\{\sigma_{2k_i+1}^{(n_i)}<\infty\}}\Big(\int_{\sigma_{2k_i+1}^{(n_i)}-T_2}^{\sigma_{2k_i+1}^{(n_i)}}\sigma^2(s)\D s\Big)^{\frac{1}{2}}\Big]\\
    &\le&r_uT_2+8\sigma_u^2T_2^{\frac{1}{2}}\\
    &<&\infty.
\end{eqnarray*}
This is a contradiction, since $\ln n_i\RA\infty$ as $i\RA\infty$. Therefore, Case 1 is not true.

If Case 2 arises, then by (\ref{T1.8}) and (\ref{TC2}) we have
\begin{eqnarray}
P\{\Omega_{2k_j+1}^{(n_j)}(T_2)\}&=&P(\Omega_{n_j})-P\big(\{\sigma_{2k_j+1}^{(n_j)}-\tau_{k_j}^{(n_j)}< T_2\}\cap\Omega_{n_j}\big)\nonumber\\
&>&\delta_0-\eta_j\quad\mbox{for all }j=1,2,\cdots.\label{P}
\end{eqnarray}
We can choose an integer $q$ such that $T_2\in[q\gamma,(q+1)\gamma)$, then from (\ref{T1.1}), (\ref{T1.4}) and (\ref{T1.11}) for any $j=1,2,\cdots$ we have
\begin{eqnarray*}
    0&<&U_j:=E\Big[I_{\Omega_{k_j}^{n_j}(T_2)}\Big(\ln x_{n_j}(\sigma_{2k_j+1}^{(n_j)})-\ln x_{n_j}(\sigma_{2k_j+1}^{(n_j)}-T_2)\Big)\Big]\\
    &\le&E\Big[I_{\Omega_{k_j}^{n_j}(T_2)}\Big(\int_{\sigma_{2k_j+1}^{(n_j)}-T_2}^{\sigma_{2k_j+1}^{(n_j)}-T_2+q\gamma}
    +\int_{\sigma_{2k_j+1}^{(n_j)}-T_2+q\gamma}^{\sigma_{2k_j+1}^{(n_j)}}\Big)(r(s)-\frac{1}{2}\sigma^2(s)-a(s)L)\D s\Big]\\
    &&+E\Big[I_{\{\sigma_{2k_j+1}^{(n_j)}<\infty\}\cap\Omega_0}\int_{\sigma_{2k_j+1}^{(n_j)}-T_2}^{\sigma_{2k_j+1}^{(n_j)}}
    \sigma(s)\D B(s)\Big]\\
    &&-E\Big[\big(I_{\{\sigma_{2k_j+1}^{(n_j)}<\infty\}\cap\Omega_0}-I_{\Omega_{k_j}^{n_j}(T_2)}\big)\int_{\sigma_{2k_j+1}^{(n_j)}-T_2}^{\sigma_{2k_j+1}^{(n_j)}}\sigma(s)\D B(s)\Big]\\
    &\le&(-q\mu+\alpha_2\gamma) P\{\Omega_{k_j}^{n_j}(T_2)\}\\
    &&+\Big|E\big[\big(I_{\{\sigma_{2k_j+1}^{(n_j)}<\infty\}\cap\Omega_0}-I_{\Omega_{n_j}}+I_{\{\sigma_{2k_j+1}^{(n_j)}
    -\tau_{k_j}^{(n_j)}<T_2\}\cap\Omega_{n_j}}\big)\int_{\sigma_{2k_j+1}^{(n_j)}-T_2}^{\sigma_{2k_j+1}^{(n_j)}}
    \sigma(s)\D B(s)\big]\Big|.
    \end{eqnarray*}
     Since $T_2>\frac{\gamma}{\mu}(\mu+\alpha_2\gamma+1)$ and  $\rho<\delta_0^2/(9\sigma_u^2T_2)$,  we follows from (\ref{T1.10}), (\ref{TC2}), (\ref{P}) and the H$\rm\ddot{o}$lder inequality have
    \begin{eqnarray*}
    U_j&\le&(-\frac{T_2\mu}{\gamma}+\mu+\alpha_2\gamma)(\delta_0-\eta_j)+\Big(P\big\{\big(\{\sigma_{2k_j+1}^{(n_j)}<\infty\}\cap\Omega_0\big)\backslash
    \Omega_{n_j}\big\}\\
    &&+P\big\{\{\sigma_{2k_j+1}^{(n_j)}
    -\tau_{k_j}^{(n_j)}<T_2\}\cap\Omega_{n_j}\big\}\Big)^\frac{1}{2}\Big(E\big[I_{\{\sigma_{2k_j+1}^{(n_j)}<\infty\}}\int_{\sigma_{2k_j+1}^{(n_j)}-T_2}^
    {\sigma_{2k_j+1}^{(n_j)}}\sigma(s)\D B(s)\big]^2\Big)^\frac{1}{2}\\
    &\le&(-\frac{T_2\mu}{\gamma}+\mu+\alpha_2\gamma)(\delta_0-\eta_j)+(\rho+\eta_j)^\frac{1}{2}E\Big[I_{\{\sigma_{2k_j+1}^{(n_j)}<\infty\}}\int_{\sigma_{2k_j+1}^{(n_j)}-T_2}^
    {\sigma_{2k_j+1}^{(n_j)}}\sigma^2(s)\D s\Big]^\frac{1}{2}\\
    &\le&-(\delta_0-\eta_j)+(\rho+\eta_j)^\frac{1}{2}\sigma_uT_2^\frac{1}{2}.
\end{eqnarray*}
Since $\eta_j\RA0$ as $j\RA\infty$, there is a $J$ such that for all $j\ge J$
$$
    \eta_j\le\min\{\rho,\frac{\delta_0}{3}\}.
$$
Hence, we obtain
$$
    U_j\le-\frac{\delta_0}{3}(2-\sqrt{2})\quad\mbox{for all }j\ge J,
$$
which contradict with $U_j>0$ for all $j=1,2,\cdots$. Therefore, Case 2 is not true and our claim (\ref{CL2}) hold.

By (\ref{CL2}), for any $T>\frac{\gamma}{\mu}(\mu+\alpha_2\gamma+1)$, there are positive constants $N>2$ and $L_n$ such that
\begin{equation}
    P\{\Omega_k^n(T)\}=P\{\Omega_n\}\quad\mbox{for all }n\ge N,\;k\ge L_n.\label{T1.12}
\end{equation}
By (\ref{T1.9}) and (\ref{T1.12}), for any positive constant $\rho<\delta_0/(4\sigma_u^2T)$ and $n\ge N$ there exist a $K'_n\ge L_n$ such that
\begin{equation}
    P\left\{\left(\{\sigma_k^{(n)}<\infty\}\cap\Omega_0\right)\backslash\Omega_k^n(T)\right\}<\rho\quad\mbox{for all }k\ge K'_n.\label{T1.10'}
\end{equation}
Choosing an integer $p>0$ such that $T\in[p\gamma,(p+1)\gamma)$, for any $n>N$ and $k>K'_n$ from (\ref{T1.4}), (\ref{T1.1}) and (\ref{T1.11}) we have
\begin{eqnarray*}
    0&<&U_k^n:=E\Big[I_{\Omega_k^n(T)}\Big(\ln x_n(\sigma_{2k+1}^{(n)})-\ln x_n(\sigma_{2k+1}^{(n)}-T)\Big)\Big]\\
    &\le&E\Big[I_{\Omega_k^n(T)}\Big(\int_{\sigma_{2k+1}^{(n)}-T}^{\sigma_{2k+1}^{(n)}-T+p\gamma}
    +\int_{\sigma_{2k+1}^{(n)}-T+p\gamma}^{\sigma_{2k+1}^{(n)}}\Big)(r(s)-\frac{1}{2}\sigma^2(s)-a(s)L)\D s\Big]\\
    &&+E\Big[I_{\{\sigma_{2k+1}^{(n)}<\infty\}\cap\Omega_0}\int_{\sigma_{2k+1}^{(n)}-T}^{\sigma_{2k+1}^{(n)}}
    \sigma(s)\D B(s)\Big]\\
    &&-E\Big[\big(I_{\{\sigma_{2k+1}^{(n)}<\infty\}\cap\Omega_0}
    -I_{\Omega_k^n(T)}\big)\int_{\sigma_{2k+1}^{(n)}-T}^{\sigma_{2k+1}^{(n)}}\sigma(s)\D B(s)\Big]\\
    &\le&(-p\mu+\alpha_2\gamma) P\{\Omega_k^n(T)\}+\Big|E\Big[I_{\big(\{\sigma_{2k+1}^{(n)}<\infty\}\cap\Omega_0\big)
    \backslash{\Omega_k^n(T)}}\int_{\sigma_{2k+1}^{(n)}-T}^{\sigma_{2k+1}^{(n)}}
    \sigma(s)\D B(s)\Big]\Big|\\
    &\le&(-\frac{T\mu}{\gamma}+\mu+\alpha_2\gamma)P\{\Omega_k^n(T)\}\\
    &&+P\Big\{\big(\{\sigma_{2k+1}^{(n)}<\infty\}\cap\Omega_0\big)
    \backslash{\Omega_k^n(T)}\Big\}^\frac{1}{2}\Big(E\big[I_{\{\sigma_{2k+1}^{(n)}<\infty\}}\int_{\sigma_{2k+1}^{(n)}-T}^
    {\sigma_{2k+1}^{(n)}}\sigma(s)\D B(s)\big]^2\Big)^\frac{1}{2}\\
    &=&(-\frac{T\mu}{\gamma}+\mu+\alpha_2\gamma)P\{\Omega_k^n(T)\}\\
    &&+P\Big\{\big(\{\sigma_{2k+1}^{(n)}<\infty\}\cap\Omega_0\big)
    \backslash{\Omega_k^n(T)}\Big\}^\frac{1}{2}E\Big[I_{\{\sigma_{2k+1}^{(n)}<\infty\}}\int_{\sigma_{2k+1}^{(n)}-T}^
    {\sigma_{2k+1}^{(n)}}\sigma^2(s)\D s\Big]^\frac{1}{2}.
    \end{eqnarray*}
    Since $T>\frac{\gamma}{\mu}(\mu+\alpha_2\gamma+1)$ and  $\rho<\delta_0^2/(4\sigma_u^2T)$, from (\ref{T1.8}), (\ref{T1.12}) and (\ref{T1.10'}) we have
    \begin{equation*}
        U_k^n\le-\delta_0+\rho^\frac{1}{2}\sigma_uT^\frac{1}{2}\le-\frac{\delta_0}{2},
    \end{equation*}
    which is a contradiction. Therefore, the conclusion of the Theorem \ref{T1} is hold. This complete the proof.\qed
\end{theorem}

\begin{remark}
\rm By (i) of Lemma \ref{L2}, we can obtain that
$$
    \LI P\{x(t)\le M\}\ge P\{\LS x(t)\le M\}.
$$
Therefore, the result of Theorem \ref{T1} is more general than the stochastically ultimately bounded.
\end{remark}
\begin{theorem}\label{T2}
    Suppose ($H_1$) and ($H_2$) hold.
    Then for any $\varepsilon\in(0,1)$ there is a positive constant $m=m(\varepsilon)$ such that for any initial value $x_0\in R_+$ the solution obeys
    $$
        P\{\LI x(t)\ge m\}\ge 1-\varepsilon.
    $$
    \rm\textbf{Proof.} From the assumption ($\rm H_2$), there are positive constants $T_0$, $l$ and $\mu$ such that
\begin{equation}
    \int_t^{t+\lambda}\big(r(s)-\frac{1}{2}\sigma^2(s)-a(s)l\big)\D s>\mu\quad\mbox{for all }t\ge T_0.\label{T2.1}
\end{equation}
Then we can choose a positive constant $\delta_0\ll1$ such that
\begin{equation}
    \int_t^{t+\lambda}\big(r(s)-\frac{1}{2}\sigma^2(s)-a(s)(l+\delta_0)\big)\D s>\frac{\mu}{2}\quad\mbox{for all }t\ge T_0.\label{T2.1'}
\end{equation}
Let
$$
\beta_1=\sup_{s\ge0}|r(s)-\frac{1}{2}\sigma^2(s)-a(s)(l+\delta_0)|\quad\mbox{and}\quad\beta_2=\sup_{s\ge0}|r(s)-\frac{1}{2}\sigma^2(s)-a(s)l|.
$$
Firstly, we will prove
\begin{equation}
    \LS x(t,x_0,\omega)>l\quad\mbox{for all }x_0\in R_+\mbox{ and }\omega\in\Omega_0.\label{CL2.1}
\end{equation}
Otherwise, there exist a $x_0\in R_+$ and $\omega_0\in\Omega_0$ such that
$$
    \LS x(t,x_0,\omega_0)\le l.
$$
Hence, for above $\delta_0$ there is a $T_1\ge T_0$ such that
\begin{equation}
    x(t,x_0,\omega_0)\le l+\delta_0\quad\mbox{ for all }t\ge T_1.\label{T2.5}
\end{equation}
For any $t\ge T_1$ we can choose a positive integer $p$ such that $t\in[T_1+p\lambda,T_1+(p+1)\lambda)$, then by (\ref{T1.4}), (\ref{T2.1'}) and (\ref{T2.5}) we get
\begin{eqnarray}
    \ln x(t)-\ln x(T_1)&\ge&\Big(\int_{T_1}^{T_1+p\lambda}+\int_{T_1+p\lambda}^t\Big)\Big(r(s)-\frac{1}{2}\sigma^2(s)-a(s)(l+\delta_0)\Big)\D s\nonumber\\
    &&+M(t)-M(T_1)\nonumber\\
    &\ge&\frac{1}{2}p\mu-\lambda\beta_1+M(t)-M(T_1)\nonumber\\
    &\ge&\frac{\mu}{2\lambda}(t-T_1-\lambda)-\lambda\beta_1+M(t)-M(T_1).\label{T2.6}
\end{eqnarray}
Dividing $t$ on the both sides and then letting $t\RA\infty$, it finally follows from (\ref{P1}) that
$$
    \LI\frac{\ln x(t)}{t}\ge\frac{\mu}{2\lambda}.
$$
This implies $\lim_{t\RA\infty}x(t)=\infty$, which is a contradiction. Therefore, (\ref{CL2.1}) hold.

In the following, we will prove this theorem is true. Otherwise, there is a positive constant $\varepsilon_0\in(0,1)$ such that for any positive integer $n\ge2$ we have a $x_n\in R_+$ satisfying the property that
$$
    P\{\LI x_n(t)\ge \frac{l}{n}\}\le1-\varepsilon_0,
$$
where $x_n(t)=x(t,x_n)$. Consequently,
\begin{equation}
    P\{\LI x_n(t)<\frac{l}{n}\}\ge \varepsilon_0.\label{T2.7}
\end{equation}
Let us now define a sequence of stopping time for any $n\ge2$,
\begin{eqnarray*}
    \sigma_1^{(n)}&=&\inf\{t\ge T_0:x_n(t)\le \frac{l}{n}\},\quad \sigma_{2k}^{(n)}=\inf\{t\ge \sigma_{2k-1}^{(n)}:x_n(t)\ge l\}\\
    \sigma_{2k+1}^{(n)}&=&\inf\{t\ge \sigma_{2k}^{(n)}:x_n(t)\le \frac{l}{n}\},\quad\mbox{for }k=1,2,\cdots.
\end{eqnarray*}
Let
$$\Omega_n=\Big(\{\LI x_n(t)<\frac{l}{n}\}\cup\big(\{\LI x_n(t)=\frac{l}{n}\}\cap\big(\bigcap_{k=1}^{\infty}\{\sigma_k^{(n)}<\infty\}\big)\big)\Big)\cap\Omega_0$$ for $n\ge2.$
 Then from (\ref{T2.7}) we obtain
\begin{equation}
    P(\Omega_n)\ge\varepsilon_0.\label{T2.8}
\end{equation}
Note from (\ref{CL2.1}) and the definition of $\Omega_n$ that for any $n=2,3,\cdots$
$$
    \sigma_k^{(n)}<\infty\quad\mbox{for any }k\ge 1\mbox{ whenever }\omega\in\Omega_n,
$$
$\{\sigma_k^{(n)}<\infty\}$ is non-increasing with $k$ and
\begin{equation}
    \Omega_n=\Big(\bigcap_{k=1}^{\infty}\{\sigma_k^{(n)}<\infty\}\Big)\cap\Omega_0.\label{T2.9}
\end{equation}
Now for any $n\ge2$ we can define a time sequence $\{\tau_k^{(n)}\}$, where
$$
    \tau_k^{(n)}=\sup\{t:x_n(t)\ge l\mbox{ and }\sigma_{2k}^{(n)}\le t<\sigma_{2k+1}^{(n)}\},\quad k=1,2,\cdots.
$$
By (\ref{CL2.1}) and the definition of $\Omega_n$, we can find that for any $\omega\in\Omega_n$
\begin{equation}
    \begin{array}{c}
        \displaystyle\frac{l}{n}<x_n(t)<l\quad\mbox{for all }t\in(\tau_k^{(n)},\sigma_{2k+1}^{(n)}),\vspace{1mm}\\
       \displaystyle x_n(\tau_k^{(n)})=l \mbox{ and }x_n(\sigma_{2k+1}^{(n)})=\frac{l}{n}\quad\mbox{for all }n\ge2,\;k=1,2,\cdots.
    \end{array}\label{T2.11}
\end{equation}
We claim that for any positive constant $T>\frac{\lambda}{\mu}(\mu+\beta_2\lambda+1)$, we have an integer $N$, for any $n\ge N$ there is a  $H_n$ such that
\begin{equation}
    P\big\{\{\sigma_{2k+1}^{(n)}-\tau_k^{(n)}\ge T\}\cap\Omega_n\big\}=P(\Omega_n)\quad\mbox{for all } k\ge H_n.\label{CL2.2}
\end{equation}
If the claim is not true, there is a positive constant $T_2>\frac{\lambda}{\mu}(\mu+\beta_2\lambda+1)$, for any $i=1,2,\cdots$ there exist a positive constat $n_i>i$ such that for any $h=1,2,\cdots$ we have a positive constant $k_i^h\ge h$ satisfying
\begin{equation}
    P\big\{\{\sigma_{2k_i^h+1}^{(n_i)}-\tau_{k_i^h}^{(n_i)}\ge T_2\}\cap\Omega_{n_i}\big\}<P(\Omega_{n_i}).\label{TTT}
\end{equation}
By (\ref{T2.9}), for any positive constant $\nu^\frac{1}{2}<\varepsilon_0/(3\sigma_uT_2^\frac{1}{2})$ and $n\ge2$ there exist a $K_n=K_n(\nu)$ such that
\begin{equation}
    P\Big\{\big(\{\sigma_k^{(n)}<\infty\}\cap\Omega_0\big)\backslash\Omega_n\Big\}<\nu\quad\mbox{for all }k\ge K_n.\label{T2.10}
\end{equation}
For $i=1,2,\cdots$ and $h\ge K_{n_i}$ then there is a $k_i^h\ge h$ satisfying (\ref{TTT}) and (\ref{T2.10}).
Consequently, it follows from (\ref{TTT}) we have
$$
    P\big\{\{\sigma_{2k_i^h+1}^{(n_i)}-\tau_{k_i^h}^{(n_i)}<T_2\}\cap\Omega_{n_i}\big\}>0.
$$
Obviously, there are two cases:\\
Case 1. There is a positive constant $\eta_0$ such that
$$
    P\big\{\{\sigma_{2k_i^h+1}^{(n_i)}-\tau_{k_i^h}^{(n_i)}<T_2\}\cap\Omega_{n_i}\big\}>\eta_0\quad\mbox{for all }i=1,2,\cdots\mbox{ and }h\ge K_{n_i};
$$
Case 2. There are subsequences $\{n_{i_j}\}\subset\{n_i\}$ and $\{k_{i_j}^{h_j}:j=1,2,\cdots,h_j\ge K_{n_{i_j}}\}\subset\{k_i^h:i=1,2,\cdots,h\ge K_{n_i}\}$, a sequence $\{\eta_j\}$ with $\eta_j\RA 0$ as $j\RA\infty$ and $\eta_j<\varepsilon_0$ for all $j=1,2, \cdots$ such that
\begin{equation}
    0<P\big\{\{\sigma_{2k_{i_j}+1}^{(n_{i_j})}-\tau_{k_{i_j}}^{(n_{i_j})}<T_2\}\cap\Omega_{n_{i_j}}\big\}<\eta_j\quad\mbox{for all }j=1,2,\cdots,\label{TC2.2}
\end{equation}
where $k_{i_j}:=k_{i_j}^{h_j}$. Here $n_{i_j}$ and $k_{i_j}^{h_j}$ may not be increasing with $j$.

In the following, we will prove neither Case 1 nor Case 2 is true. If Case 1 arises, from Theorem \ref{T1} for above $\eta_0$ there is a positive constant $M=M(\eta_0)$ such that
$$
    P\{\LS x(t,x_0)\le M\}\ge1-\frac{\eta_0}{2}\quad\mbox{for any }x_0\in R_+.
$$
Therefore, we can obtain that
$$
     P\big\{\{\sigma_{2k_i^h+1}^{(n_i)}-\tau_{k_i^h}^{(n_i)}<T_2\}\cap\Omega_{n_i}\cap\{\LS x_{n_i}(t)\le M\}\big\}>\frac{\eta_0}{2},
$$
for all $i=1,2,\cdots$ and $h\ge K_{n_i}$.
We denote
$$\Delta_{k_i^h}^{(n_i)}=\{\sigma_{2k_i^h+1}^{(n_i)}-\tau_{k_i^h}^{(n_i)}<T_2\}\cap\Omega_{n_i}\cap\{\LS x_{n_i}(t)\le M\}$$
and
$$\Delta^{(n_i)}=\Omega_{n_i}\cap\{\LS x_{n_i}(t)\le M\},$$
then
$$P\{\Delta^{(n_i)}\}\ge P\{\Delta_{k_i^h}^{(n_i)}\}>\frac{\eta_0}{2}.$$
Consequently, follows from (\ref{T1.4}), (\ref{T1.11}) and the moment inequality of stochastic integrals (see \cite{M97}), we compute
\begin{eqnarray*}
    -\frac{\eta_0}{2}\ln n_i&=&-\frac{\eta_0}{2}(\ln l-\ln \frac{l}{n_i})\\
    &\ge&-E\Big[\sup_{0\le t\le T_2}I_{\Delta_{k_i^h}^{(n_i)}}\left(\ln x_{n_i}(\sigma_{2k_i^h+1}^{(n_i)}-t)-\ln x_{n_i}(\sigma_{2k_i^h+1}^{(n_i)})\right)\Big]\\
    &\ge&-E\Big[\sup_{0\le t\le T_2}I_{\Delta_{k_i^h}^{(n_i)}}\int_{\sigma_{2k_i^h+1}^{(n_i)}-t}^{\sigma_{2k_i^h+1}^{(n_i)}}\left(a(s)x(s)+\frac{1}{2}\sigma^2(s)+r(s)\right)\D s\Big]\\
    &&-E\Big[\sup_{0\le t\le T_2}I_{\Delta_{k_i^h}^{(n_i)}}\Big|\int_{\sigma_{2k_i^h+1}^{(n_i)}-t}^{\sigma_{2k_i^h+1}^{(n_i)}}\sigma(s)\D B(s)\Big|\Big]\\
    &\ge&-E\Big[I_{\Delta^{(n_i)}}\int_{\sigma_{2k_i^h+1}^{(n_i)}-T_2}^{\sigma_{2k_i^h+1}^{(n_i)}}a(s)x(s)\D s\Big]-(\frac{1}{2}\sigma_u^2+r_u)T_2\\
    &&-E\Big[I_{\{\sigma_{2k_i^h+1}^{(n_i)}<\infty\}}\sup_{0\le t\le T_2}\Big|\int_{\sigma_{2k_i^h+1}^{(n_i)}-t}^{\sigma_{2k_i^h+1}^{(n_i)}}\sigma(s)\D B(s)\Big|\Big]\\
    &\ge&-E\Big[I_{\Delta^{(n_i)}}\int_{\sigma_{2k_i^h+1}^{(n_i)}-T_2}^{\sigma_{2k_i^h+1}^{(n_i)}}a(s)x(s)\D s\Big]\\
    &&-(\frac{1}{2}\sigma_u^2+r_u)T_2-8E\Big[I_{\{\sigma_{2k_i^h+1}^{(n_i)}<\infty\}}\int_{\sigma_{2k_i^h+1}^{(n_i)}-T_2}^{\sigma_{2k_i^h+1}^{(n_i)}}\sigma^2(s)\D s\Big]^{\frac{1}{2}}\\
    &\ge&-E\Big[I_{\Delta^{(n_i)}}\int_{\sigma_{2k_i^h+1}^{(n_i)}-T_2}^{\sigma_{2k_i^h+1}^{(n_i)}}a(s)x(s)\D s\Big]-(\frac{1}{2}\sigma_u^2+r_u)T_2-8\sigma_uT_2^{\frac{1}{2}}.
\end{eqnarray*}
Since $\sigma_{2k_i^h+1}^{(n_i)}\RA \infty$ as $h\RA\infty$ for all $\omega\in\Omega_{n_i}$, from above inequality, letting $h\RA\infty$ we have
\begin{eqnarray*}
    -\frac{\eta_0}{2}\ln n_i&\ge&-\limsup_{h\RA\infty}E\Big[I_{\Delta^{(n_i)}}\int_{\sigma_{2k_i^h+1}^{(n_i)}-T_2}^{\sigma_{2k_i^h+1}^{(n_i)}}a(s)x(s)\D s\Big]\\
    &&-(\frac{1}{2}\sigma_u^2+r_u)T_2-8\sigma_uT_2^{\frac{1}{2}}\\
    &\ge&-E\Big[I_{\Delta^{(n_i)}}\limsup_{h\RA\infty}\int_{\sigma_{2k_i^h+1}^{(n_i)}-T_2}^{\sigma_{2k_i^h+1}^{(n_i)}}a(s)x(s)\D s\Big]\\
    &&-(\frac{1}{2}\sigma_u^2+r_u)T_2-8\sigma_uT_2^{\frac{1}{2}}\\
    &\ge&-E\Big[I_{\Delta^{(n_i)}}Ma_uT_2\Big]-(\frac{1}{2}\sigma_u^2+r_u)T_2-8\sigma_uT_2^{\frac{1}{2}}\\
    &\ge&-(Ma_u+\frac{1}{2}\sigma_u^2+r_u)T_2-8\sigma_uT_2^{\frac{1}{2}}.
\end{eqnarray*}
This is a contradiction, since $\ln n_i\RA\infty$ as $i\RA\infty$. Therefore, Case 1 is not true.

If Case 2 arises, then by (\ref{T2.8}) and (\ref{TC2.2}) we have
\begin{eqnarray}
&&P\left(\{\sigma_{2k_{i_j}+1}^{(n_{i_j})}-\tau_{k_{i_j}}^{(n_{i_j})}\ge T_2\}\cap\Omega_{n_{i_j}}\right)\nonumber\\
    &=&P(\Omega_{n_{i_j}})-P\left(\{\sigma_{2k_{i_j}+1}^{(n_{i_j})}-\tau_{k_{i_j}}^{(n_{i_j})}< T_2\}\cap\Omega_{n_{i_j}}\right)\nonumber\\
    &>&\varepsilon_0-\eta_j\quad\mbox{for all }j=1,2,\cdots.\label{P2}
\end{eqnarray}
We can choose an integer $q$ such that $T_2\in[q\lambda,(q+1)\lambda)$, then from (\ref{T1.4}), (\ref{T2.1}) and (\ref{T2.11}) for any $j=1,2,\cdots$ we have
\begin{eqnarray*}
    0&>&V_j:=E\Big[I_{\{\sigma_{2k_{i_j}+1}^{(n_{i_j})}-\tau_{k_{i_j}}^{(n_{i_j})}\ge T_2\}\cap\Omega_{n_{i_j}}}\Big(\ln x_{n_{i_j}}(\sigma_{2k_{i_j}+1}^{(n_{i_j})})-\ln x_{n_{i_j}}(\sigma_{2k_{i_j}+1}^{(n_{i_j})}-T_2)\Big)\Big]\\
    &\ge&E\Big[I_{\{\sigma_{2k_{i_j}+1}^{(n_{i_j})}-\tau_{k_{i_j}}^{(n_{i_j})}\ge T_2\}\cap\Omega_{n_{i_j}}}\Big(\int_{\sigma_{2k_{i_j}+1}^{(n_{i_j})}-T_2}^{\sigma_{2k_{i_j}+1}^{(n_{i_j})}-T_2+q\lambda}
    +\int_{\sigma_{2k_{i_j}+1}^{(n_{i_j})}-T_2+q\lambda}^{\sigma_{2k_{i_j}+1}^{(n_{i_j})}}\Big)\\
    &&(r(s)-\frac{1}{2}\sigma^2(s)-a(s)l)\D s\Big]+E\Big[I_{\{\sigma_{2k_{i_j}+1}^{(n_{i_j})}<\infty\}\cap\Omega_0}\int_{\sigma_{2k_{i_j}+1}^{(n_{i_j})}-T_2}^{\sigma_{2k_{i_j}+1}^{(n_{i_j})}}
    \sigma(s)\D B(s)\Big]\\
    &&-E\Big[\big(I_{\{\sigma_{2k_{i_j}+1}^{(n_{i_j})}<\infty\}\cap\Omega_0}-I_{\{\sigma_{2k_{i_j}+1}^{(n_{i_j})}-\tau_{k_{i_j}}^{(n_{i_j})}\ge T_2\}\cap\Omega_{n_{i_j}}}\big)\int_{\sigma_{2k_{i_j}+1}^{(n_{i_j})}-T_2}^{\sigma_{2k_{i_j}+1}^{(n_{i_j})}}\sigma(s)\D B(s)\Big]\\
    &\ge&(q\mu-\beta_2\lambda) P\{\{\sigma_{2k_{i_j}+1}^{(n_{i_j})}-\tau_{k_{i_j}}^{(n_{i_j})}\ge T_2\}\cap\Omega_{n_{i_j}}\}\\
    &&-\Big|E\Big[\big(I_{\{\sigma_{2k_{i_j}+1}^{(n_{i_j})}<\infty\}\cap\Omega_0}-I_{\Omega_{n_{i_j}}}+I_{\{\sigma_{2k_{i_j}+1}^{(n_{i_j})}
    -\tau_{k_{i_j}}^{(n_{i_j})}<T_2\}\cap\Omega_{n_{i_j}}}\big)\int_{\sigma_{2k_{i_j}+1}^{(n_{i_j})}-T_2}^{\sigma_{2k_{i_j}+1}^{(n_{i_j})}}
    \sigma(s)\D B(s)\Big]\Big|.
    \end{eqnarray*}
     Since $T_2>\frac{\lambda}{\mu}(\mu+\beta_2\lambda+1)$ and  $\nu^\frac{1}{2}<\varepsilon_0/(3\sigma_uT_2^\frac{1}{2})$,  we follows from (\ref{T2.10}), (\ref{TC2.2}), (\ref{P2}) and the H$\rm\ddot{o}$lder inequality have
    \begin{eqnarray*}
    V_j&\ge&(\frac{T_2\mu}{\lambda}-\mu-\beta_2\lambda)(\varepsilon_0-\eta_j)-\Big(P\big\{\{\sigma_{2k_{i_j}+1}^{(n_{i_j})}<\infty\}\cap\Omega_0\backslash
    \Omega_{n_{i_j}}\big\}\\
    &&+P\big\{\{\sigma_{2k_{i_j}+1}^{(n_{i_j})}
    -\tau_{k_{i_j}}^{(n_{i_j})}<T_2\}\cap\Omega_{n_{i_j}}\big\}\Big)^\frac{1}{2}\Big(E\big[I_{\{\sigma_{2k_{i_j}+1}^{(n_{i_j})}<\infty\}}\int_{\sigma_{2k_{i_j}+1}^{(n_{i_j})}-T_2}^
    {\sigma_{2k_{i_j}+1}^{(n_{i_j})}}\sigma(s)\D B(s)\big]^2\Big)^\frac{1}{2}\\
    &\ge&(\frac{T_2\mu}{\lambda}-\mu-\beta_2\lambda)(\varepsilon_0-\eta_j)-(\nu+\eta_j)^\frac{1}{2}E\Big[I_{\{\sigma_{2k_{i_j}+1}^{(n_{i_j})}<\infty\}}\int_{\sigma_{2k_{i_j}+1}^{(n_{i_j})}-T_2}^
    {\sigma_{2k_{i_j}+1}^{(n_{i_j})}}\sigma^2(s)\D s\Big]^\frac{1}{2}\\
    &\ge&(\varepsilon_0-\eta_j)-(\nu+\eta_j)^\frac{1}{2}\sigma_uT_2^\frac{1}{2}.
\end{eqnarray*}
Since $\eta_j\RA0$ as $j\RA\infty$, there is a $J$ such that for all $j\ge J$
$$
    \eta_j\le \min\{\nu,\frac{\varepsilon_0}{3}\}.
$$
Hence, we obtain
$$
    V_j\ge\frac{\varepsilon_0}{3}( 2-\sqrt{2})\quad\mbox{for all }j\ge J,
$$
which contradict with $V_j<0$ for all $j=1,2,\cdots$. Therefore, Case 2 is not true and our claim (\ref{CL2.2}) hold.

By (\ref{CL2.2}), for any $T>\frac{\lambda}{\mu}(\mu+\beta_2\lambda+1)$, there is a positive constant $N>2$ such that for any $n\ge N$ we have a positive integer $H_n$ satisfying
\begin{equation}
    P\{\Omega_k^n(T)\}=P\{\Omega_n\}\quad\mbox{for all }k\ge H_n,\label{T2.12}
\end{equation}
where
$$\Omega_k^n(T)=\{\sigma_{2k+1}^{(n)}-\tau_k^{(n)}\ge T\}\cap\Omega_n.$$
By (\ref{T2.9}) and (\ref{T2.12}), for any positive constants $\nu^\frac{1}{2}<\delta_0/(2\sigma_uT^\frac{1}{2})$ and $n\ge N$ there exist a $K'_n\ge H_n$ such that
\begin{equation}
    P\left\{\left(\{\sigma_k^{(n)}<\infty\}\cap\Omega_0\right)\backslash\Omega_k^n(T)\right\}<\nu\quad\mbox{for all }k\ge K'_n.\label{T2.10'}
\end{equation}
Choosing an integer $p>0$ such that $T\in[p\lambda,(p+1)\lambda)$, for any $n>N$ and $k>K'_n$ from (\ref{T1.4}), (\ref{T2.1}) and (\ref{T2.11}) we have
\begin{eqnarray*}
    0&>&V_k^{(n)}:=E\Big[I_{\Omega_k^n(T)}\Big(\ln x_n(\sigma_{2k+1}^{(n)})-\ln x_n(\sigma_{2k+1}^{(n)}-T)\Big)\Big]\\
    &\ge&E\Big[I_{\Omega_k^n(T)}\Big(\int_{\sigma_{2k+1}^{(n)}-T}^{\sigma_{2k+1}^{(n)}-T+p\lambda}
    +\int_{\sigma_{2k+1}^{(n)}-T+p\lambda}^{\sigma_{2k+1}^{(n)}}\Big)(r(s)-\frac{1}{2}\sigma^2(s)-a(s)l)\D s\Big]\\
    &&+E\Big[I_{\{\sigma_{2k+1}^{(n)}<\infty\}\cap\Omega_0}\int_{\sigma_{2k+1}^{(n)}-T}^{\sigma_{2k+1}^{(n)}}
    \sigma(s)\D B(s)\Big]\\
    &&-E\Big[\big(I_{\{\sigma_{2k+1}^{(n)}<\infty\}\cap\Omega_0}
    -I_{\Omega_k^n(T)}\big)\int_{\sigma_{2k+1}^{(n)}-T}^{\sigma_{2k+1}^{(n)}}\sigma(s)\D B(s)\Big]\\
    &\ge&(p\mu-\beta_2\lambda) P\{\Omega_k^n(T)\}-\Big|E\Big[I_{\big(\{\sigma_{2k+1}^{(n)}<\infty\}\cap\Omega_0\big)
    \backslash{\Omega_k^n(T)}}\int_{\sigma_{2k+1}^{(n)}-T}^{\sigma_{2k+1}^{(n)}}
    \sigma(s)\D B(s)\Big]\Big|\\
    &\ge&(\frac{T\mu}{\lambda}-\mu-\beta_2\lambda)P\{\Omega_k^n(T)\}\\
    &&-\Big[P\left\{\big(\{\sigma_{2k+1}^{(n)}<\infty\}\cap\Omega_0\big)
    \backslash{\Omega_k^n(T)})\right\}\Big]^\frac{1}{2}\Big(E\big[I_{\{\sigma_{2k+1}^{(n)}<\infty\}}\int_{\sigma_{2k+1}^{(n)}-T}^
    {\sigma_{2k+1}^{(n)}}\sigma(s)\D B(s)\big]^2\Big)^\frac{1}{2}\\
    &=&(\frac{T\mu}{\lambda}-\mu-\beta_2\lambda)P\{\Omega_k^n(T)\}\\
    &&-\Big[P\left\{\big(\{\sigma_{2k+1}^{(n)}<\infty\}\cap\Omega_0\big)
    \backslash{\Omega_k^n(T)})\right\}\Big]^\frac{1}{2}\Big(E\big[I_{\{\sigma_{2k+1}^{(n)}<\infty\}}\int_{\sigma_{2k+1}^{(n)}-T}^
    {\sigma_{2k+1}^{(n)}}\sigma^2(s)\D s\big]\Big)^\frac{1}{2}.
    \end{eqnarray*}
    Since $T>\frac{\lambda}{\mu}(\mu+\beta_2\lambda+1)$ and  $\nu^\frac{1}{2}<\varepsilon_0/(2\sigma_uT^\frac{1}{2})$, from (\ref{T2.8}), (\ref{T2.12}) and (\ref{T2.10'}) we have
    \begin{eqnarray*}
        V_k^{(n)}&\ge&\varepsilon_0-\nu^\frac{1}{2}\sigma_uT^\frac{1}{2}\\
        &\ge&\frac{\varepsilon_0}{2},
    \end{eqnarray*}
    which is a contradiction. Therefore, the conclusion of the Theorem \ref{T2} is hold. This complete the proof.\qed
\end{theorem}
By (ii) of Lemma \ref{L2} and Theorem \ref{T2}, we can obtain the following result immediately.
\begin{theorem}\label{T2'}
     Suppose the assumptions ($H_1$) and ($H_2$) hold.
    Then for any $\varepsilon\in(0,1)$ there is a positive constant $m=m(\varepsilon)$ such that
    $$
        \LI P\{ x(t)\ge m\}\ge 1-\varepsilon,
    $$
    for any positive solution $x(t)$ of system (\ref{S1}).
\end{theorem}
Consequently, from Theorem \ref{T1}-\ref{T2'} we have that
\begin{theorem}\label{T5}
     Suppose the assumptions ($H_1$) and ($H_2$) hold.
     Then system (\ref{S1}) is stochastically permanence. Furthermore, for any $\varepsilon\in(0,1)$ there are positive constants $m=m(\varepsilon)$ and $M=M(\varepsilon)$ such that for any initial value $x_0\in R_+$ the solution obeys
    \begin{equation}
        P\{\LI x(t)\ge m\}\ge 1-\varepsilon\mbox{ and }P\{\LS x(t)\le M\}\ge 1-\varepsilon.\label{C1.1}
    \end{equation}
\end{theorem}
\begin{remark}
    \rm In this section, we have used three theorems to complete the proof of the stochastically permanence of system (\ref{S1}) and provided a new method for studying the stochastically permanence of stochastic differential equation, which is completely different from that of \cite{JSL} and \cite{LM}. Furthermore, \cite{JSL} obtained the stochastically permanence of system (\ref{S1}) under the conditions that $r(t)$, $a(t)$ and $\sigma(t)$ are positive continuous $T$-periodic functions and $\min_{t\in[0,T]}r(t)>\max_{t\in[0,T]}\sigma^2(t)$. In \cite{LM}, the authors studied the stochastically permanence of system (\ref{S1}) with $a_l>0$ and $(r-\frac{1}{2}\sigma^2)_l>0$. In \cite{LW}, the authors show that if
    $$
        \LI (r(t)-\frac{1}{2}\sigma^2(t))>0\mbox{ and } a_l>0,
    $$
     then equation (\ref{S1}) is stochastically permanent. Obviously, the conditions ($\rm H_1$) and ($\rm H_2$) are more weaker than these, and the result (\ref{C1.1}) is more generally than the stochastically permanence of system (\ref{S1}), that is system (\ref{S1}) is stochastically permanent while (\ref{C1.1}) is not always true. Hence, Theorem \ref{T5} is more general than these and we have the following corollary.
\end{remark}
\begin{corollary}\label{R1}
    Suppose that $a_l>0$ and $\LI (r-\frac{1}{2}\sigma^2)_l>0$. Then system (\ref{S1}) is stochastically permanence. Furthermore, for any $\varepsilon\in(0,1)$ there are positive constants $m=m(\varepsilon)$ and $M=M(\varepsilon)$ such that for any initial value $x_0\in R_+$ the solution obeys
    $$
        P\{\LI x(t)\ge m\}\ge 1-\varepsilon\mbox{ and }P\{\LS x(t)\le M\}\ge 1-\varepsilon.
    $$
\end{corollary}
\begin{remark}\label{R2}
    \rm In Liu and Wang (2011), the authors obtained the strongly persistent in the mean a.s., that is
    $$
        \LI \frac{1}{t}\int_0^t x(t)\D s>0 \mbox{ a.s.}
    $$
    with the conditions
    $$
        \LI \frac{1}{t}\int_0^t(r(t)-\frac{1}{2}\sigma^2(t))\D s>0\mbox{ and } a_l>0.
    $$
    It is well known that if $f(t)$ is an almost periodic function, then the condition
    $$
        \lim_{t\RA+\infty}\frac{1}{t}\int_0^tf(s)\D s>(<)0
    $$
    is equivalent to the condition that there is a positive constant $\lambda$ such that
    $$
        \liminf_{t\RA+\infty}\int_t^{t+\lambda}f(s)\D s>(<)0.
    $$
    On the other hand, stochastically permanence implies strong persistence in the mean. Therefore, for the particular case if $r(t)$, $a(t)$ and $\sigma(t)$ are continuous almost periodic functions, the conditions in Theorem \ref{T5} are much weaker, while the results are much better than these and we have following corollary.
\end{remark}
\begin{corollary}\label{AC1}
    Suppose $r(t)$, $a(t)$ and $\sigma(t)$ are continuous almost periodic functions,
    $$
        \lim_{t\RA\infty}\frac{1}{t}\int_0^t(r(s)-\frac{1}{2}\sigma^2(s))\D s>0\mbox{ and }\lim_{t\RA\infty}\frac{1}{t}\int_0^ta(s)\D s>0.
    $$
    Then system (\ref{S1}) is stochastically permanence. Furthermore, for any $\varepsilon\in(0,1)$ there are positive constants $m=m(\varepsilon)$ and $M=M(\varepsilon)$ such that for any initial value $x_0\in R_+$ the solution obeys
    $$
        P\{\LI x(t)\ge m\}\ge 1-\varepsilon\mbox{ and }P\{\LS x(t)\le M\}\ge 1-\varepsilon.
    $$
\end{corollary}

\begin{remark}\label{m1}
    \rm The Lemma 1 of \cite{TL} studied the permanence of logistic equation, they obtained that if there are positive constants $\lambda$ and $\gamma$ such that
    $$
        \LI \int_t^{t+\lambda}r(t)\D s>0\mbox{ and }\LI \int_t^{t+\gamma}a(s)\D s>0.
    $$
    Then the deterministic system (\ref{S2}) is permanence. Hence, from Theorem \ref{T5} we can find that under certain conditions the original non-autonomous equation (\ref{S2}) and the associated stochastic equation (\ref{S1}) have similar dynamic behave such as ultimately bounded and permanent. In other words, we show that under certain conditions the noise will not spoil these nice properties.
\end{remark}

\section{Extinction}
In this section, we study the extinction of system (\ref{S1}).
\begin{theorem}\label{T3}
    Suppose ($\rm H_1$) and ($\rm H_2$) hold.
    Then for any positive solution $x(t)$ of equation (\ref{S1}) we have
    $$
        \lim_{t\RA\infty}x(t)=0\;a.s.
    $$
    \rm\textbf{Proof.} By ($\rm H_1$), there are positive constants $\delta_0$ and $T_0$ such that
    $$
        \int_t^{t+\gamma}a(s)\D s\ge \delta_0\quad\mbox{for all }t\ge T_0.
    $$
    Letting $Q=[\frac{\gamma}{\lambda}]+1$. From ($\rm H_3$), for any number $\varepsilon\in(0,\frac{\delta_0}{2})$ there is a positive constant $T_1=T_1(\varepsilon)\ge T_0$ such that
    $$
        \int_t^{t+Q\lambda}\big(r(s)-\frac{1}{2}\sigma^2(s)\big)\D s\le\varepsilon^2\quad\mbox{for all }t\ge T_1.
    $$
    Hence, we have
    \begin{eqnarray}
         &&\int_t^{t+Q\lambda}\big(r(s)-\frac{1}{2}\sigma^2(s)-\varepsilon a(s)\big)\D s\nonumber\\
         &\le&\varepsilon^2-\varepsilon\int_t^{t+\gamma} a(s)\D s\nonumber\\
         &\le&(\varepsilon-\delta_0)\varepsilon\nonumber\\
         &\le&-\frac{\delta_0}{2}\varepsilon\quad\mbox{for all }t\ge T_1.\label{T3.2}
    \end{eqnarray}
    Then by using a similar argument as the discussion about (\ref{CL1}) in Theorem \ref{T1} we can obtain  $x_0\in R_+$
    $$
        \LI x(t,x_0,\omega)<\varepsilon\quad\mbox{for all }x_0\in R_+\mbox{ and }\omega\in\Omega_0.
    $$
    Therefore, by the arbitrariness of $\varepsilon$ we have
    $$
        \LI x(t,x_0,\omega)=0\quad\mbox{for all }x_0\in R_+\mbox{ and }\omega\in\Omega_0,
    $$
    i.e.
    \begin{equation}
        P\{\LI x(t,x_0)=0\}=1\quad\mbox{for all }x_0\in R_+.\label{T3.3}
    \end{equation}
    If this theorem is not true, then there is a $x_0\in R_+$ and a positive constant $\varepsilon_0\in(0,\frac{\delta_0}{2})$ such that inequality (\ref{T3.2}) hold with $\varepsilon=\varepsilon_0$ and $T_1=T_1(\varepsilon_0)$ and
    \begin{equation}
        P\{\LS x(t)>\varepsilon_0\}\ge\varepsilon_0,\label{T3.4}
    \end{equation}
    where $x(t)$ is the solution of equation (\ref{S1}) with initial value $x(0)=x_0$. Similar argument as Theorem \ref{T2} we also can define a sequence of stopping time for any $n\ge2$,
\begin{eqnarray*}
    \sigma_1^{(n)}&=&\inf\{t\ge T_1:x(t)\le \frac{\varepsilon_0}{n}\},\quad \sigma_{2k}^{(n)}=\inf\{t\ge \sigma_{2k-1}^{(n)}:x(t)\ge \varepsilon_0\}\\
    \sigma_{2k+1}^{(n)}&=&\inf\{t\ge \sigma_{2k}^{(n)}:x(t)\le \frac{\varepsilon_0}{n}\},\quad\mbox{for }k=1,2,\cdots.
\end{eqnarray*}
Let
$$\Omega_n=\bigg(\{\LS x(t)>\varepsilon_0\}\cup\Big(\{\LS x(t)=\varepsilon_0\}\cap\big(\bigcap_{k=1}^{\infty}\{\sigma_k^{(n)}<\infty\}\big)\Big)\bigg)\cap\Omega_0$$ for $n\ge2.$
 It follows from (\ref{T3.3}) and (\ref{T3.4}) that
\begin{equation}
    P(\Omega_n)\ge\varepsilon_0.\label{T3.5}
\end{equation}
Note from (\ref{T3.3}) and the definition of $\Omega_n$ that for any $n=2,3,\cdots$
$$
    \sigma_k^{(n)}<\infty\quad\mbox{for any }k\ge 1\mbox{ whenever }\omega\in\Omega_n,
$$
$\{\sigma_k^{(n)}<\infty\}$ is non-increasing with $k$ and
\begin{equation*}
    \Omega_n=\Big(\bigcap_{k=1}^{\infty}\{\sigma_k^{(n)}<\infty\}\Big)\cap\Omega_0.
\end{equation*}
For any $n\ge2$ we can also define a time sequence $\{\tau_k^{(n)}\}$, where
$$
    \tau_k^{(n)}=\sup\{t:x_n(t)\ge \varepsilon_0\mbox{ and }\sigma_{2k}^{(n)}\le t<\sigma_{2k+1}^{(n)}\},\quad k=1,2,\cdots.
$$
By the definition of $\Omega_n$, we can find that for any $\omega\in\Omega_n$
\begin{equation*}
    \begin{array}{c}
        \displaystyle\frac{\varepsilon_0}{n}<x_n(t)<\varepsilon_0\quad\mbox{for all }t\in(\tau_k^{(n)},\sigma_{2k+1}^{(n)}),\vspace{1 mm}\\
       \displaystyle x_n(\tau_k^{(n)})=\varepsilon_0 \mbox{ and }x_n(\sigma_{2k+1}^{(n)})=\frac{\varepsilon_0}{n}\quad\mbox{for all }n\ge2,\;k=1,2,\cdots.
    \end{array}
\end{equation*}
The following discussion is rather similar with Theorem \ref{T2} and we can obtain (\ref{T3.4}) is not arising. Thus, the proof is completed.\qed
\end{theorem}
\begin{corollary}
    Suppose $(r-\frac{1}{2}\sigma^2)_u\le0$ and $a_l>0$. Then
    $$
        \lim_{t\RA\infty}x(t)=0\mbox{ a.s.}
    $$
\end{corollary}
\begin{remark}
    \rm In this theorem, we introduce a new method for studying the extinction of stochastic equation. In Li and Mao (2009) and Liu and Wang (2011), the authors obtained the extinction of system (\ref{S1}) under the conditions
    $$
        \LS\frac{1}{t}\int_0^t(r(s)-\frac{1}{2}\sigma^2(s))\D s<0\mbox{ and }a_l>0.
    $$
    But for the case
    $$
       \LS\frac{1}{t}\int_0^t(r(s)-\frac{1}{2}\sigma^2(s))\D s=0,
    $$
    Liu and Wang (2011) only obtained the system (\ref{S1}) is non-persistence in the mean a.s., that is
    $$
        \lim_{t\RA\infty}\frac{1}{t}\int_0^tx(s)\D s=0\mbox{ a.s.}
    $$
    It is well known that for any continuous function $f(t)$ if there is a positive constant $\lambda$ such that
    $$
        \LS\int_t^{t+\lambda}f(s)\D s\le0,
    $$
    then we can obtain
    $$
        \LS\frac{1}{t}\int_0^tf(s)\D s\le0.
    $$
    Hence, for some special case if
    $$
       \LS\frac{1}{t}\int_0^t(r(s)-\frac{1}{2}\sigma^2(s))\D s=0,
    $$
    we can obtain the extinction of system (\ref{S1}).
    Therefore, Theorem \ref{T3} is a new result and improved \cite{LW} extensively and we have following corollary.
\end{remark}
\begin{corollary}\label{AC2}
    Suppose $r(t)$, $a(t)$ and $\sigma(t)$ are almost periodic functions and satisfying the following conditions
    $$
        \lim_{t\RA\infty}\frac{1}{t}\int_0^t(r(s)-\frac{1}{2}\sigma^2(s))\D s<0\mbox{ and }\lim_{t\RA\infty}\frac{1}{t}\int_0^t a(s)\D s>0.
    $$
    Then
    $$
        \lim_{t\RA\infty}x(t)=0\mbox{ a.s.,}
    $$
    for any positive solution $x(t)$ of system (\ref{S1}).
\end{corollary}
\begin{remark}
\rm From Remark \ref{m1} and Theorem \ref{T3}, we find that if the noise is sufficiently large, the stochastic equation (\ref{S1}) will be extinction with probability one, although the solution to the original equation (\ref{S2}) may be persistent. In other words, the theorem reveals the important fact that environmental noise may make the population extinct.
\end{remark}
\begin{remark}
   \rm Suppose $r(t)$, $a(t)$ and $\sigma(t)$ are almost periodic functions. From Remark \ref{R2}, Theorem \ref{T5} and Theorem \ref{T3}, we can find that there is a threshold between permanence and extinction of system (\ref{S1}). That is, under the condition
   $$
    \lim_{t\RA\infty}\frac{1}{t}\int_0^t a(s)\D s>0,
   $$
   if
   $$
   \lim_{t\RA\infty}\frac{1}{t}\int_0^t(r(s)-\frac{1}{2}\sigma^2(s))\D s>0,
   $$
   system (\ref{S1}) is permanence, and if
   $$
   \lim_{t\RA\infty}\frac{1}{t}\int_0^t(r(s)-\frac{1}{2}\sigma^2(s))\D s<0,
   $$
   we have system (\ref{S1}) is extinction.
\end{remark}
\section{Global attractivity}
On the global asymptotical stability of positive solutions for system (\ref{S1}), we have the following
result.
\begin{theorem}\label{T4}
    Suppose ($\rm H_1$) and ($\rm H_2$) hold.
    Then system (\ref{S1}) is globally attractive, that is
    $$
        \lim_{t\RA\infty}|x(t)-y(t)|=0\;a.s.
    $$
    for any two positive solutions $x(t)$ and $y(t)$ of equation (\ref{S1}).\\
    \rm\textbf{Proof.} Let $x(t)$ and $y(t)$ be any two solutions of equation (\ref{S1}) with initial values $x(0)$, $y(0)\in R_+$. From Theorem \ref{T1} and Theorem \ref{T2}, we have that for any positive constant $\varepsilon<\frac{1}{4}$ there are positive constants $M=M(\varepsilon)$ and $m=m(\varepsilon)$ such that
    $$
        P\{\LS x(t)\le M\}\ge 1-\varepsilon,\quad P\{\LI x(t)\ge m\}\ge 1-\varepsilon,
    $$
    $$
        P\{\LS y(t)\le M\}\ge 1-\varepsilon\mbox{ and }P\{\LI y(t)\ge m\}\ge 1-\varepsilon.
    $$
    Consequently,
    we obtain
    \begin{equation}
        P\{\Lambda\}\ge 1-4\varepsilon>0,\label{T4.1}
    \end{equation}
    where
    $$\Lambda:=\{m\le \LI x(t)\le\LS x(t)\le M,m\le \LI y(t)\le\LS y(t)\le M\}.$$
    By the It$\rm\hat{o}$ formula, we have
    $$
        \D \ln x(t)=\big(r(t)-\frac{1}{2}\sigma^2(t)-a(t)x(t)\big)\D t+\sigma(t)\D B(t)
    $$
    and
    $$
        \D \ln x(t)=\big(r(t)-\frac{1}{2}\sigma^2(t)-a(t)x(t)\big)\D t+\sigma(t)\D B(t).
    $$
    Then
    \begin{equation}
        \D(\ln x(t)-\ln y(t))=-a(t)(x(t)-y(t))\D t.\label{T4.2}
    \end{equation}
    Choose the Lyapunov function as follows
    $$
        V(t)=|\ln x(t)-\ln y(t)|\quad\mbox{for all }t\ge 0.
    $$
    By calculating the right differential $\D^+V(t)$ of $V(t)$ along ordinary differential equation (\ref{T4.2}) leads to
    \begin{eqnarray}
        \D^+V(t)&=&\mbox{sgn}(x(t)-y(t))\D(\ln x(t)-\ln y(t))\nonumber\\
        &=&-a(t)|x(t)-y(t)|\D t.\label{T4.3}
    \end{eqnarray}
    For any $\omega\in\Lambda$ we can denote $l(\omega)=\inf_{t\ge0}\{x(t,\omega),y(t,\omega)\}$ and $L(\omega)=\sup_{t\ge0}\{x(t,\omega),y(t,\omega)\}$. By (\ref{T4.1}), we obviously have
    $$
        0<l(\omega)\le L(\omega)<\infty\quad\mbox{for all }\omega\in\Lambda.
    $$
    From (\ref{T4.3})and the famous Cauchy's Mean-Value Theorem we have that
    \begin{eqnarray*}
        \D^+V(t,\omega)&=&-a(t)\xi(t,\omega)|\ln x(t,\omega)-\ln y(t,\omega)|\D t\\
        &\le&-a(t)l(\omega)V(t,\omega)\D t\quad\mbox{for all }\omega\in\Lambda,
    \end{eqnarray*}
    where $\xi(t,\omega)$ is between $x(t,\omega)$ and $y(t,\omega)$. Consequently,
    $$
        V(t,\omega)\le V(0)\exp\Big(-\int_0^ta(s)l(\omega)\D s\Big)\quad\mbox{for all }\omega\in\Lambda.
    $$
    Since $\int_0^\infty a(s)\D s=+\infty$, we have $V(t,\omega)\RA0$ as $t\RA\infty$ for all $\omega\in\Lambda$. Combine with (\ref{T4.1}) we get
    $$
        P\{\lim_{t\RA\infty} V(t)=0\}\ge P\{\Lambda\}\ge 1-4\varepsilon.
    $$
    Then by the arbitrariness of $\varepsilon$ we easily know that
    $$
        P\{\lim_{t\RA\infty} V(t)=0\}=1.
    $$
    Therefore,
    $$
        \lim_{t\RA\infty}|x(t)-y(t)|=0\;a.s.
    $$
This completes the proof.\qed
\end{theorem}
\begin{remark}
    \rm In Teng and Li (2000), the authors studied the global attractivity of the deterministic system (\ref{S2}) under the same conditions of the permanence of the system. Hence, from Theorem \ref{T4} we can find that under certain conditions the stochastic equation (\ref{S1}) keep the global attractivity property.
\end{remark}
\begin{remark}
    \rm In Li and Mao (2009), the authors obtained the global atrractivity of system (\ref{S2}) with the condition $a_l>0$. Obviously, our result is different with their. From Theorem \ref{T3} and Theorem \ref{T4}, we can obtain that under assumption ($\rm H_1$) if assumption ($\rm H_2$) or ($\rm H_3$) hold, then system (\ref{S1}) is global attractivity. But we have not discussed the case if for any positive constant $\lambda$
    $$
        \LI\int_t^{t+\lambda} (r(s)-\frac{1}{2}\sigma^2(s))\D s<0\mbox{ and }\LS\int_t^{t+\lambda} (r(s)-\frac{1}{2}\sigma^2(s))\D s>0.
    $$
    If $r(t)$, $a(t)$ and $\sigma(t)$ are almost periodic functions, then we have not discussed the case
    $$
        \lim_{t\RA\infty}\frac{1}{t}\int_0^t (r(s)-\frac{1}{2}\sigma^2(s))\D s=0.
    $$
   Therefore, there is an interesting open question that if the assumption ($\rm H_1$) hold whether the system (\ref{S2}) is global attractivity.

\end{remark}

\section{Examples}
In this section we will give four examples to illustrate the conclusions obtained in the above sections.\\
\textbf{Example 1.} Let
$$
    r(t)=\sin t+\frac{2}{3},\;a(t)=\cos t+1\mbox{ and }\sigma(t)=\sqrt{\cos t+1}.
$$
Obviously,
$$(r-\frac{1}{2}\sigma^2)_l=\frac{1}{6}-\frac{\sqrt{5}}{2}<0\mbox{ and }a_l=0.$$
Hence, the conclusions of \cite{JSL,LM,LW} can not be used. On the other hand, choosing $\lambda=\gamma=2\pi$ then we have
$$
    \int_t^{t+\lambda}(r(s)-\frac{1}{2}\sigma^2(s))\D s=\frac{\pi}{3}\mbox{ and }\int_t^{t+\gamma}a(s)\D s=2\pi.
$$
Therefore, the assumption ($\rm H_1$) and ($\rm H_2$) are satisfied, and by the Theorem \ref{T5} we can obtain that system (\ref{S1}) is permanence (see Fig.1).\\
\textbf{Example 2.} Let
$$
    r(t)=\sin t+\frac{1}{2},\;a(t)=\cos t+1\mbox{ and }\sigma(t)=\sqrt{\cos t+1}.
$$
Consequently,
$$
     (r-\frac{1}{2}\sigma^2)_l=-\frac{\sqrt{5}}{2}<0,\;(r-\frac{1}{2}\sigma^2)_u=\frac{\sqrt{5}}{2}>0,\vspace{-2 mm}
$$
$$
     \lim_{t\RA\infty}\frac{1}{t}\int_0^t(r(s)-\frac{1}{2}\sigma^2(s))\D s=0\mbox{ and }a_l=0.
$$
Hence, all the criteria for the extinction in \cite{JSL,LM,LW} will become invalid. But if we choose $\lambda=\gamma=2\pi$ we can find that
$$
    \int_t^{t+\lambda}(r(s)-\frac{1}{2}\sigma^2(s))\D s=0\mbox{ and }\int_t^{t+\gamma}a(s)\D s=2\pi\mbox{ for all }t\ge0.
$$
Follows from Theorem \ref{T3} we obtain that the system (\ref{S1}) is extinction (see Fig.2).
\begin{center}
\includegraphics[width=.50\textwidth,height=60mm]{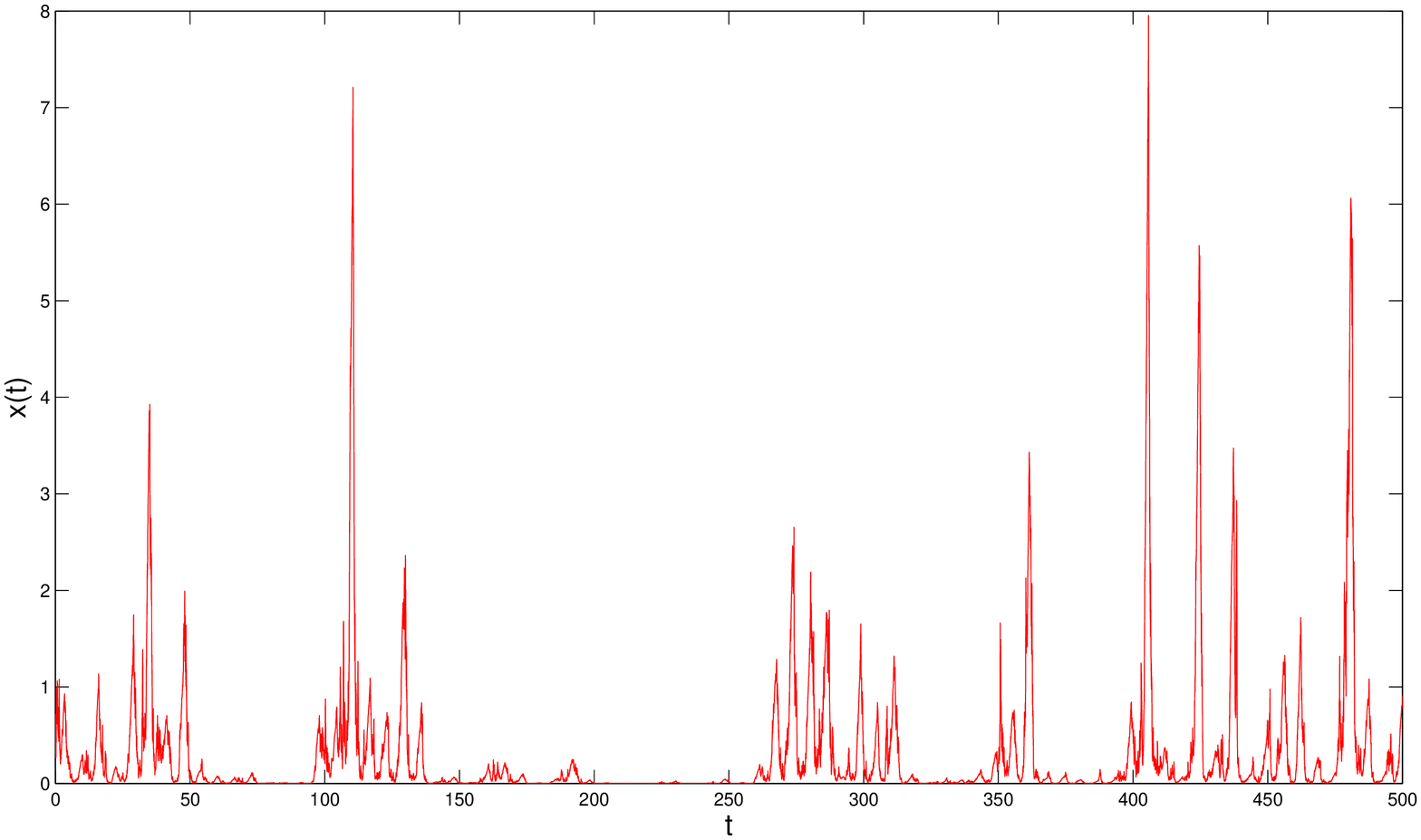}
    \put(-190,-10){\mbox{\small Fig.1. Solutions of system (\ref{S1}) for }}
    \put(-160,-30){\mbox{\small $x(0)=0.5$ in Example 1.}}
     \includegraphics[width=.50\textwidth,height=60mm]{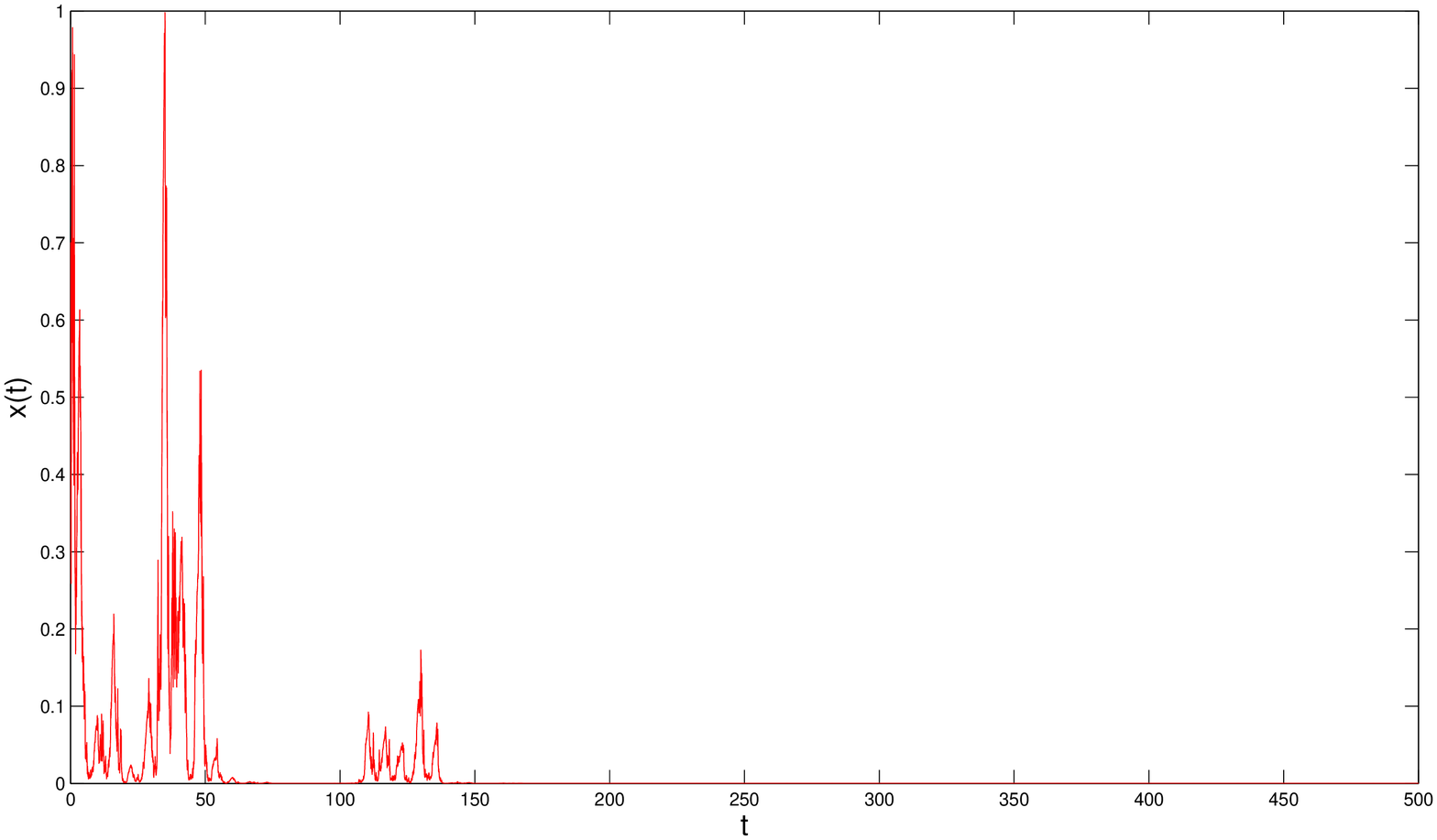}
    \put(-190,-10){\mbox{\small Fig.2. Solutions of system (\ref{S1}) for }}
    \put(-160,-30){\mbox{\small $x(0)=0.5$ in Example 2.}}
\end{center}

For some almost periodic case it is difficult to find the $\lambda$ and $\gamma$, we can use the Corollary \ref{AC1} and Corollary \ref{AC2}.\\
\textbf{Example 3.} Let
$$
    r(t)=\sin\sqrt{2}t+\cos\sqrt{3}t+\frac{2}{3},\;a(t)=\sin\sqrt{6}t+\cos\sqrt{2}t+2\mbox{ and }\sigma(t)=\sqrt{\cos t+1}.
$$
Obviously,
$$
     (r-\frac{1}{2}\sigma^2)_l=-\frac{7}{3}<0\mbox{ and }a_l=0.
$$
And it is difficult to find the $\lambda$ and $\gamma$, but we have that
$$\lim_{t\RA\infty}\frac{1}{t}\int_0^t(r(s)-\frac{1}{2}\sigma^2(s))\D s=\frac{1}{6}>0$$
and
$$\lim_{t\RA\infty}\frac{1}{t}\int_0^ta(s)\D s=2>0.$$
Therefore, follows from Corollary \ref{AC1} we obtain that the system (\ref{S1}) is permanence (see Fig.3).\\
\textbf{Example 4.}
Let
$$
    r(t)=\sin\sqrt{2}t+\cos\sqrt{3}t+\frac{1}{3},\;a(t)=\sin\sqrt{6}t+\cos\sqrt{2}t+2\mbox{ and }\sigma(t)=\sqrt{\cos t+1}.
$$
It is difficult to find the positive constants $\lambda$ and $\gamma$, but we find that
$$\lim_{t\RA\infty}\frac{1}{t}\int_0^t(r(s)-\frac{1}{2}\sigma^2(s))\D s=-\frac{1}{6}<0$$
and
$$\lim_{t\RA\infty}\frac{1}{t}\int_0^ta(s)\D s=2>0.$$
Therefore, by Corollary \ref{AC2} we have the system (\ref{S1}) is extinction (see Fig.4).
\begin{center}
\includegraphics[width=.50\textwidth,height=60mm]{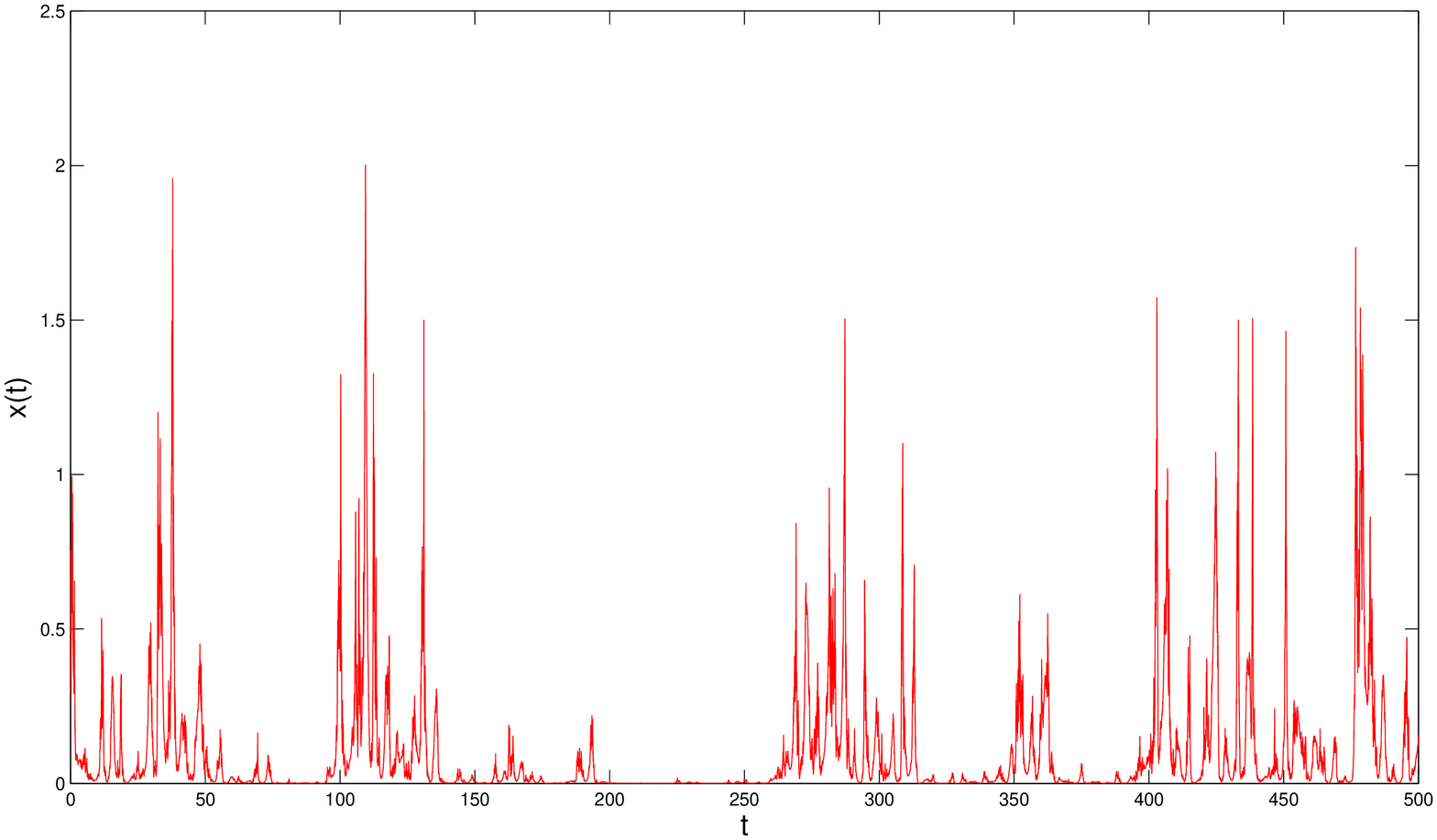}
    \put(-190,-10){\mbox{\small Fig.3. Solutions of system (\ref{S1}) for }}
    \put(-160,-30){\mbox{\small $x(0)=0.5$ in Example 3.}}
     \includegraphics[width=.50\textwidth,height=60mm]{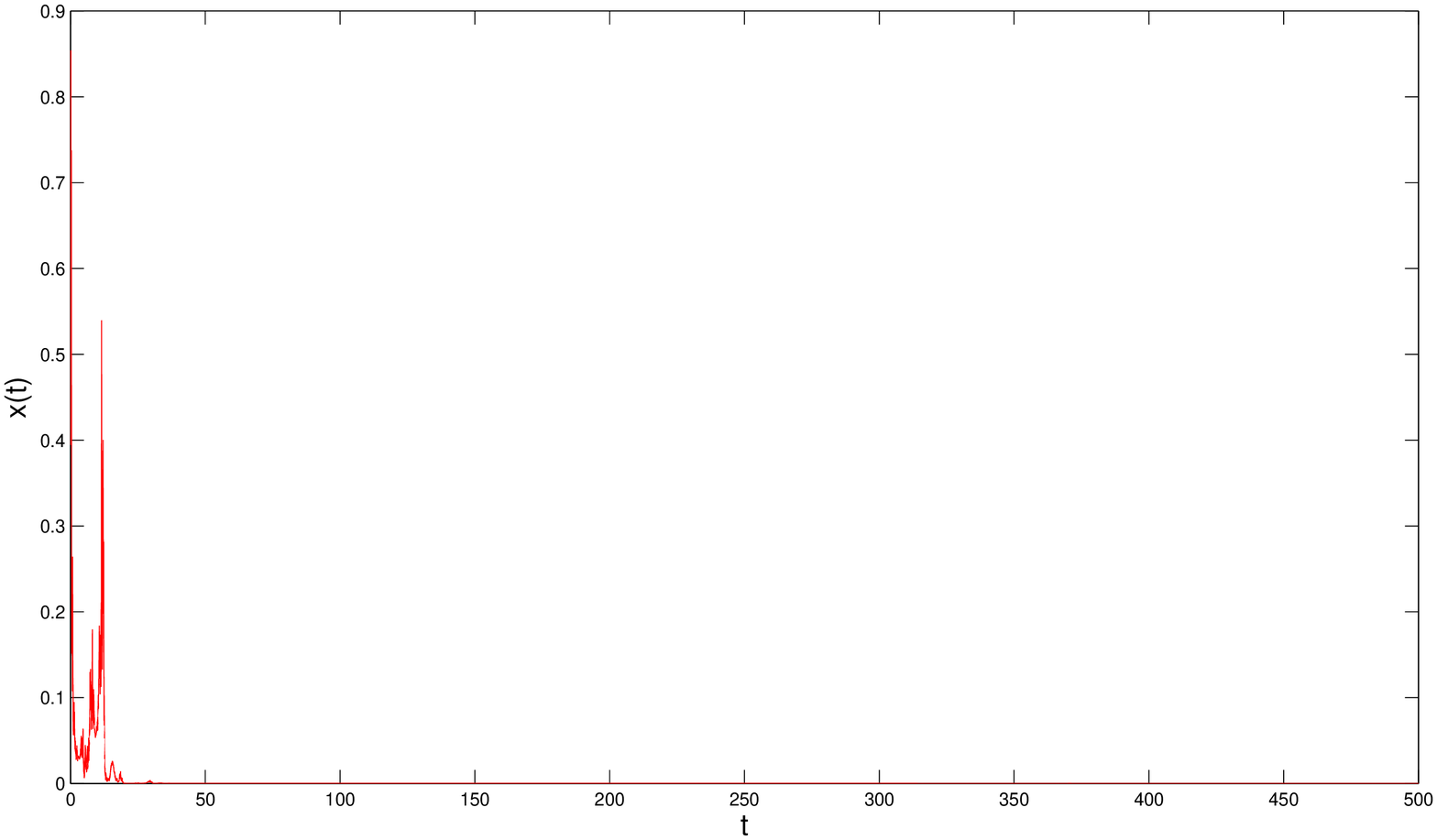}
    \put(-190,-10){\mbox{\small Fig.4. Solutions of system (\ref{S1}) for }}
    \put(-160,-30){\mbox{\small $x(0)=0.5$ in Example 4.}}
\end{center}
\bibliographystyle{plain}

\bibliography{mybio}		

\end{document}